\documentclass[preprint,12pt]{elsarticle}

\usepackage{amssymb}

\usepackage{graphicx}
\usepackage{multicol,multirow}
\usepackage{amsmath,amssymb,amsfonts}
\usepackage{mathrsfs}
\usepackage{amsthm}
\usepackage{makecell}
\usepackage{rotating}
\usepackage{appendix}
\usepackage[numbers]{natbib}
\usepackage{ifpdf}
\usepackage[T1]{fontenc}
\usepackage{newtxtext}
\usepackage{newtxmath}
\usepackage{textcomp}
\usepackage{xcolor}
\usepackage{lscape}
\usepackage{lipsum}
\usepackage[colorlinks,allcolors=blue]{hyperref}
\usepackage{caption}
\usepackage{subcaption}

\theoremstyle{definition}

\numberwithin{equation}{section}

\journal{Mathematical Biosciences}

\begin{document}

\begin{frontmatter}



\title{Clusters of African countries based on the social contacts and associated socioeconomic indicators relevant to the spread of the epidemic.}


\author[inst1]{Evans Kiptoo Korir}

\affiliation[inst1]{organization={Bolyai Institute},
            addressline={University of Szeged}, 
            postcode={6720}, 
            state={Szeged},
            country={Hungary}}
\author[inst1]{Zsolt Vizi}

\affiliation[inst1]{organization={Bolyai Institute},
            addressline={National Laboratory for Health Security}, 
            postcode={6720}, 
            state={Szeged},
            country={Hungary}}

\begin{abstract}
{
\textbf{Introduction}. 

It is well known that social contact patterns differ from country to country. This variation coincides with significant socioeconomic heterogeneity that complicates the design of effective non-pharmaceutical interventions. This study examined how socioeconomic heterogeneity in selected African countries might be factored in to explain better social contact mix patterns between countries.

\textbf{Methods}. 
We used a standardized contact matrix for 32 African countries, estimated in \cite{PK}. We scaled the matrices using an epidemic model from \cite{RB}. We also analyzed aggregated data from the World Bank country website. The data includes 28 variables; social, economic, environmental, institutional, governance, health and well-being, education, gender inequality, and other development-related indicators describing countries. Principal components analysis was used to visualize socioeconomic similarities between countries and identify the indicators for maximum variation. The $(2D)^2 PCA$ approach was used to reduce the dimension of the synthetic contact matrices for each country to avoid the dimensionality curse. Agglomerative hierarchical clustering was then used to identify groups of countries with similar social patterns, taking into account the country's socioeconomic performance.

\textbf{Results}. 
Our model yielded four meaningful clusters, each with a few distinguishing features. Social contacts varied between groups but were generally similar within each set. The country's socioeconomic performance influenced the clusters.

\textbf{Conclusions}. Our results suggest that integrating socioeconomic factors into social contacts can better explain infectious disease transmission dynamics and that similar interventions can be implemented in countries within the cluster.}

\end{abstract}




\end{frontmatter}



\section{Introduction}

The outbreak of a pandemic affects the lives of many people around the world. Influenza, measles, and the new COVID-19 are some reported infectious diseases. Transmission is mainly via droplet infection, therefore dependent on the frequency of contact between susceptible and infected persons \cite{PK}. Because of the transmissions, public health interventions, including non-pharmaceutical interventions (NPIs) such as school closures, lockdowns, and social distancing, have been implemented to change contact patterns and thereby lower the peak of the epidemic. Mathematical models, specifically compartment models, where the population is compartmentalized for the epidemic evolution simulations, have been used to simulate the outcomes and effectiveness of constraints \cite{ZA}. Since transmission events cannot be observed and measured, various assumptions are made to estimate the age-specific transmission parameters, for example, the hypothesis of homogeneous mixing and the social contact hypothesis \cite{WA}.

In particular, it is known that the frequency of contact between susceptible and infected individuals varies depending on factors such as age, gender, behavior, and the context of setting such as work, school, home, or community. Interaction rates between age groups in the environment are modeled using age-stratified models. We can infer from the model which age groups are most susceptible to the disease, as well as those with high rates of transmission. Age assortative (people tend to mingle with individuals in their age group) has been observed in several studies to estimate contact matrices; see \cite{EZ, KK, KR, MX, MH, PK, RB}. According to \cite{MH}, based on 7290 participants in 8 European countries, the POLYMOD study offers the first comprehensive quantitative approach to assessing social contact patterns. \cite{PC} estimated synthetic contact matrices for 152 nations in 2017 using the POLYMOD study, demographic information, surveys, and World Bank databases as references. Using the most recent data, the authors increased the number of nations on the list to 177. A large-scale contact database, the BBC Pandemic Project \cite{KK}, reports 36,155 participants in the UK with a total of 378,559 contacts and estimates age-specific population contact matrices in different settings. \cite{FA} has also prepared such estimates for 26 European countries. Several authors have also estimated contact matrices for each country, including \cite{AE, GG, HT, IT, KT, KG, LD, MD, RJ}.

Regional variations exist between nations, with some countries having outbreaks first and later spreading to other nations. A rapid public health response is required in each country to minimize transmission by implementing non-pharmaceutical and health measures such as vaccination. Because it is crucial to decide whether individuals comply with planned interventions, evaluating such interventions requires quantifying social relationships and considering people's livelihoods in different countries. Although introducing NPIs has proven beneficial in containing the epidemic, it has harmed the economy and people's livelihoods. The initiative has a symmetric impact on global gross domestic product (GDP) and affects the most vulnerable groups in society \cite{BR}.

In a pandemic, the government tries to find a balance between the economy's state and the disease's spread. The social distancing effort is a domestic public intervention that includes closing malls and schools, banning gatherings, and working from home. It has a more significant impact on the epidemic and the economic situation. The younger and the older population are more likely to be affected by these measures \cite{OS, SE}. Due to this adoption, nations that depend on manufacturing and agriculture are on the losing side, which disrupts supply chains and leads to high unemployment rates and lost income. Developed countries are less affected by these measures than developing countries because of digitization and their ability to use other suitable alternatives to products and services and reduce consumption, which requires interaction with suppliers and producers. Governments must implement NPIs immediately to achieve better economic outcomes, and lower mortality rates \cite{DE}. Low mobility theoretically reduces social interactions and lowers cases but may affect people's livelihood \cite{MA}. In this way, social contacts in countries with a low standard of living can accelerate the spread of the disease. As a result, additional economic indicators are needed to improve human health, promote self-sustaining momentum, and highlight the voluntary component critical to minimizing transmission and meeting NPIs.

Although the average number of contacts varies from country to country, the mixing patterns are strikingly consistent according to \cite{MH}. Policymakers monitor and evaluate the actions initiated and taken by actively intervening countries during a pandemic \cite{MX}. Before a pandemic, it is essential to know which nations to keep an eye on. We can group countries based on their connected social contacts and socioeconomic indicators important to the NPI strategy to comprehend the ones with comparable patterns and subsequently build similar NPI initiatives. We use dimensionality reduction and unsupervised machine learning approaches (principal component analysis, $(2D)^2 PCA$, and agglomerative hierarchical clustering with complete linkage method) since we have high-dimensional data.   To offset the spread of the pandemic in the country, a range of socioeconomic data was taken into account, including gross domestic product (GDP), population, mortality, life expectancy, and other indicators. A transmission model makes the contact matrices in the selected countries comparable.

\section{Methods}

\subsection{Data}

In 2017, \cite{PC} created age-specific synthetic contact matrices for 152 countries based on the POLYMOD study and data capturing household, school, and workplace information. With the most recent and updated data, the authors revised the matrices and countries to 177 (about 97 percent of the world's countries). The contact matrices are reproducible and cover rural and urban environments \cite{PK}. We have used these constructed social contact matrices but only targeting all African countries. The contact matrices were presented in four settings (home, school, work, and other) distributed across 16 age groups represented by 5-year age groups beginning with 0-4, 5-9, 10-14, $\dots$, 75+. Of the data, 32 countries in Africa have complete data and are therefore used for the analysis (see appendix for the complete list). 

\subsection{Socioeconomic indicators}
We also considered socioeconomic indicators for the selected countries from the World Bank database \cite{Z}. The data capture vital indicators such as social, economic, environmental, institutional, governance, health and well-being, education, and gender inequality. Twenty-eight indicators were considered, distributed across the key sectors, and consistent for the countries (see appendix).

\subsection{Contact Matrices}
As indicated in Fig. \ref{fig: contacts}, we established country setting-specific contact matrices for use at home, school, work, and other locations labeled $\mu_{C}^{H}, \mu_{C}^{S}, \mu_{C}^{W},$ and $\mu_{C}^{O}$ respectively. 

To arrive at the \textit{complete social contact matrix}, we consider a situation of no intervention by taking the unweighted sum of contacts made in the four settings:
\begin{equation*}
     \mu_{C} =\mu_{C}^{H} + \mu_{C}^{S} + \mu_{C}^{W} + \mu_{C}^{O}   
\end{equation*}
where $\mu_{C}$ is the complete social contact matrix for the country $C \in \{1, . . . , 32\}$.
The contacts are reciprocal; therefore, the sum of contacts from age group $i$ to $j$ must match the sum of contacts from age group $j$ to age group $i$ (see, for example, \cite{EZ, KK, KR, MX, MH, PK, RB, WH}). To maintain consistency and reciprocity, we have changed the contact matrices for each country.
We can define $W_C^{(i)}$ and $W_C^{(j)}$ as the entire population in age classes $i$ and $j$ aligned to the country $C$, respectively, as the sample has a different population than the country as a whole.

\begin{equation}
     \mu ^{(i,j)}_{C} W_C^{(i)} = \mu ^{(j, i)}_{C} W_C^{(j)} 
    \label{eq2}
\end{equation}
where  $\mu ^{(i, j)}_{C}$ is the average number of contacts a person in age group $i$ has with people in age group $j$ on a given day in the country $C$. 

We make sure that the contact matrices satisfy Equation \ref{eq2} since we wish to compare the contact matrices across various African countries. As a result, we symmetrically align the country matrices with guaranteeing contacts reciprocity.

\begin{equation*}
     \mu_{C}^{(i, j)} = \frac{1}{2 W_C^{(i)}}\left(\mathcal{\mu} ^{(i, j)}_{C} W_C^{(i)} + \mathcal{\mu} ^{(j, i)}_{C} W_C^{(j)}\right)  
\end{equation*}
where $\mu_{C}$ now exhibits strong diagonal contacts and is symmetric if we multiply by the respective population size (assortative nature). Before the clustering technique, the resulting contact matrices will be used in the subsequent dimensionality reduction stage.

\subsection{Aggregation method to transform age contact matrix}

We outline the process of generating contact matrices in the desired age format using known contact mix data and country demographics.
If $M_C$ and $N_C$ are the original contact matrix of size $16 \times 16$ and the population vector aligned with the age bins in the contact matrix for a country $C$, 
$\mu_C$ and $W_C$ referring to the new contact matrix and population vector, then
$$\mu_C^{(i,j)} = \frac{1} {W_C^{(i)}}\sum_{(m, n) \in \mathcal{I}_i \times \mathcal{I}_j} M_C^{(m,n)} N_C^{(m)}$$
where $\mathcal{I}_i$ denotes the indices of the age bins from $N_C$, which were merged into the $i^{th}$ age bin of $W_C$, e.g. $\mathcal{I}_2 = \{2, 3\},$ since $W_2$ corresponds to 5-14, thus $2^{nd}$ and $3^{rd}$ age bins (associated with age groups 5-9 and 10-14) of $N_C$ are needed for this aggregation.
The new age structure $i \in \{1, . . . , 6\}$ was chosen to provide higher resolution for all age groups in different countries.
As a demonstration of the age structure of the population of the countries $\left(W_C^{(i)}\right)$, we have used the Kenyan population in Table\ref{table:1}.

\begin{table}[h!]
\caption{Details of the age groups used in transmission model from Fig. \ref{img: diag} and  clustering of the countries}
\centering
\addtolength{\tabcolsep}{-3pt}
\begin{tabular}{c c c c c c c c} 
\Xhline{3\arrayrulewidth} 
Age class (years)& 0–4 & 5-14 & 15-19 & 20-24 & 25-64 & 65+\\ 
Age class index (i)& 1 & 2 & 3 & 4 & 5 & 6\\   
Age class $(W_C^{(i)})$  & 7044364 & 13705769 & 6010656	& 5236593 & 20424690 & 1349228 \\
population \\
\Xhline{3\arrayrulewidth}
\end{tabular}
\label{table:1}
\end{table}

\subsection{Age dependent parameters aggregation}
We elucidate the procedure for aggregating age-specific parameters of the model. We first informed the model \ref{img: diag} with the age-specific parameters from existing studies \cite{Bra, EZ, Kim, RB, ZE}.
If $N_C$ and $P$ are the original population vector and 
model parameters aligned with the age bins in the contact matrix of size $16 \times 16$ for a country $C$
and $\rho$ refers to the new model parameters then
$$\rho ^{(i)} = \frac{\sum_{m \in \mathcal{I}_i} N_{C}^{(m)} \cdot P^{(m)}} {\sum_{m \in \mathcal{I}_i}{N_C^{(m)}}}.$$
The original population vector $N_C$ and the new 6-age group population $W_C$ used in aggregating the parameters are from the Kenyan population.
All aggregated age-dependent parameters $\rho$ are listed in Table \ref{table:params}.

\begin{table}[h!]
\caption{Age-dependent model parameter estimates, same for all countries and aggregated using the Kenyan population.}
\centering
\addtolength{\tabcolsep}{-3pt}
\begin{tabular}{c c c c c c c c} 
\Xhline{3\arrayrulewidth}
\textbf{Probability} & & \textbf{Age class}\\
\cline{3-8}
 & & \textbf{0–4} & \textbf{5-14} & \textbf{15-19} & \textbf{20-24} & \textbf{25-64} & \textbf{65+}\\ 
 \Xhline{3\arrayrulewidth}
Asymptomatic course & $\theta^{i}$ & 0.95 & 0.95 & 0.9 & 0.85 & 0.85 & 0.8  \\   
Fatal outcome & $\eta^{i}$ & 0.00195 & 0.00195 & 0.0057 & 0.0057 & 0.0057 & 0.038 \\   
Intensive care  & $\zeta^{i}$ &
0.1766 & 0.1766 & 0.1766 & 0.3192 & 0.3192 & 0.3731
\\ 
(given hospitalization) \\
Hospitalization  & $h^{i}$ & 0.0205 & 0.0205 & 0.0205 & 0.1755 & 0.1755 & 0.253\\
(or intensive care) \\
\Xhline{3\arrayrulewidth}
\end{tabular}

\label{table:params}
\end{table}

\begin{figure}{}
     \centering
     \begin{subfigure}[b]{0.19\textwidth}
         \centering
         \includegraphics[width=\textwidth]{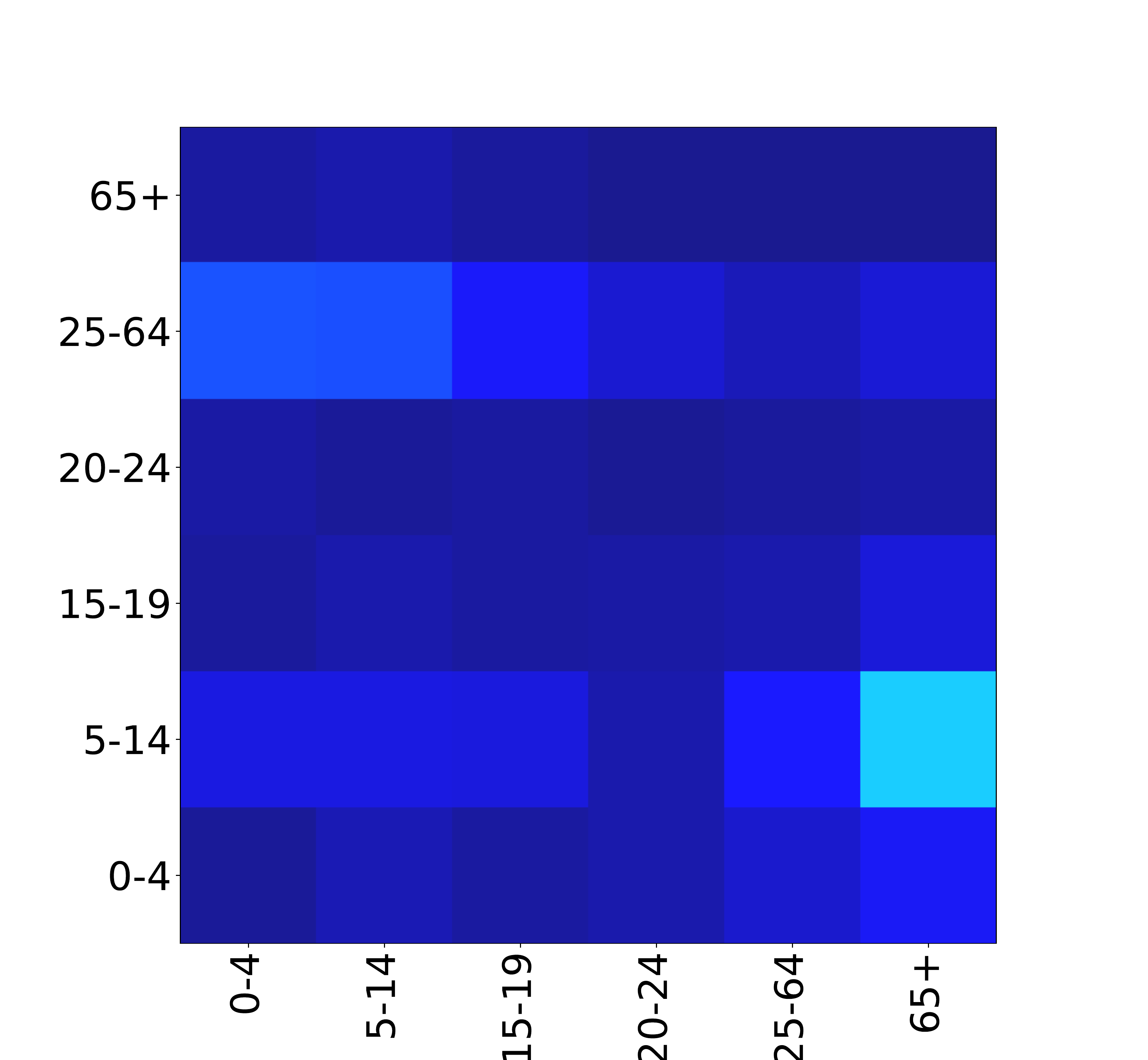}
         \caption{}
         \label{fig: home}
     \end{subfigure}
     \hfill
     \begin{subfigure}[b]{0.19\textwidth}
         \centering
         \includegraphics[width=\textwidth]{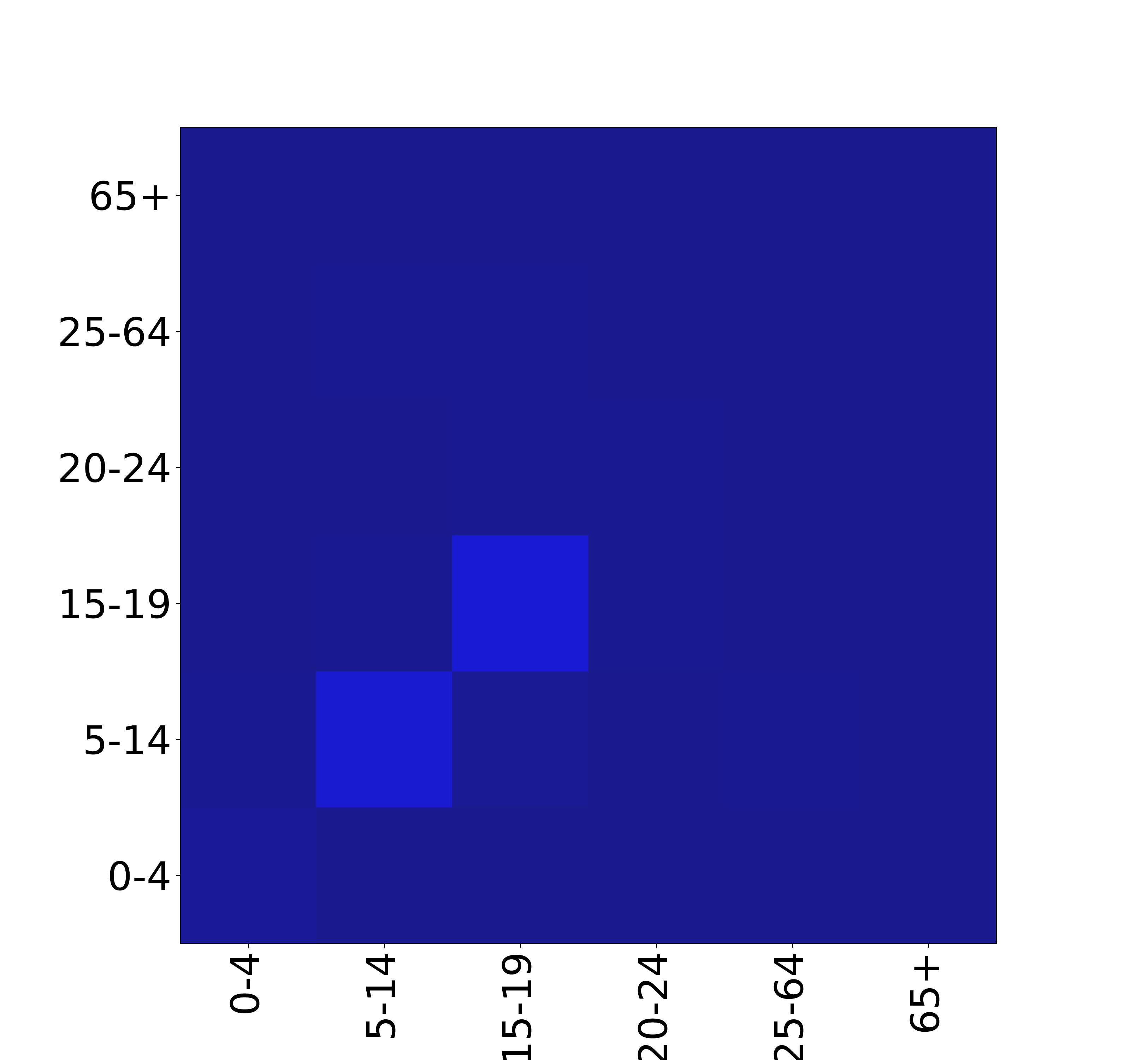}
         \caption{}
         \label{fig: school}
     \end{subfigure}
     \hfill
     \begin{subfigure}[b]{0.19\textwidth}
         \centering
         \includegraphics[width=\textwidth]{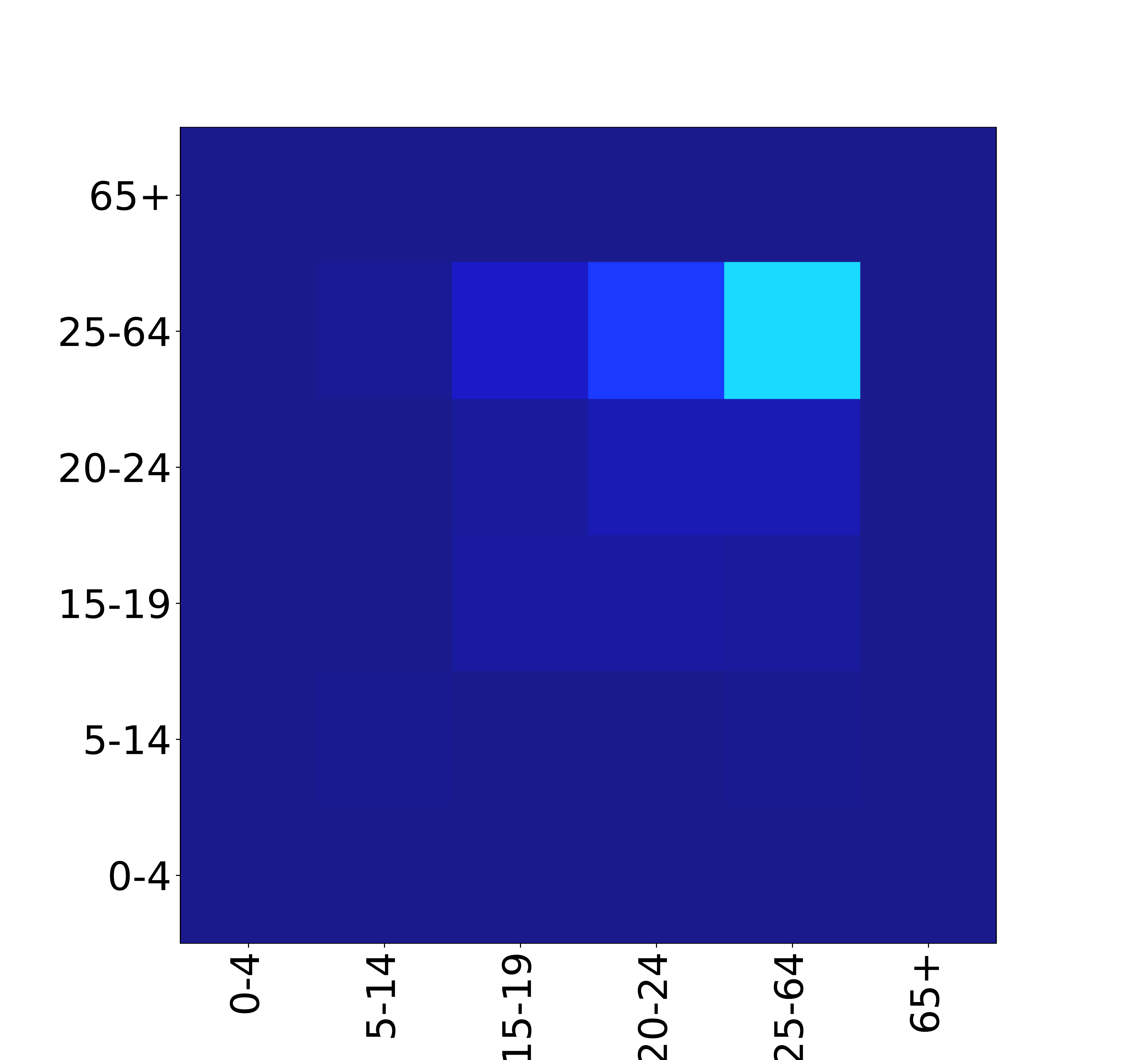}
         \caption{}
         \label{fig:work}
     \end{subfigure}
     \begin{subfigure}[b]{0.19\textwidth}
         \centering
         \includegraphics[width=\textwidth]{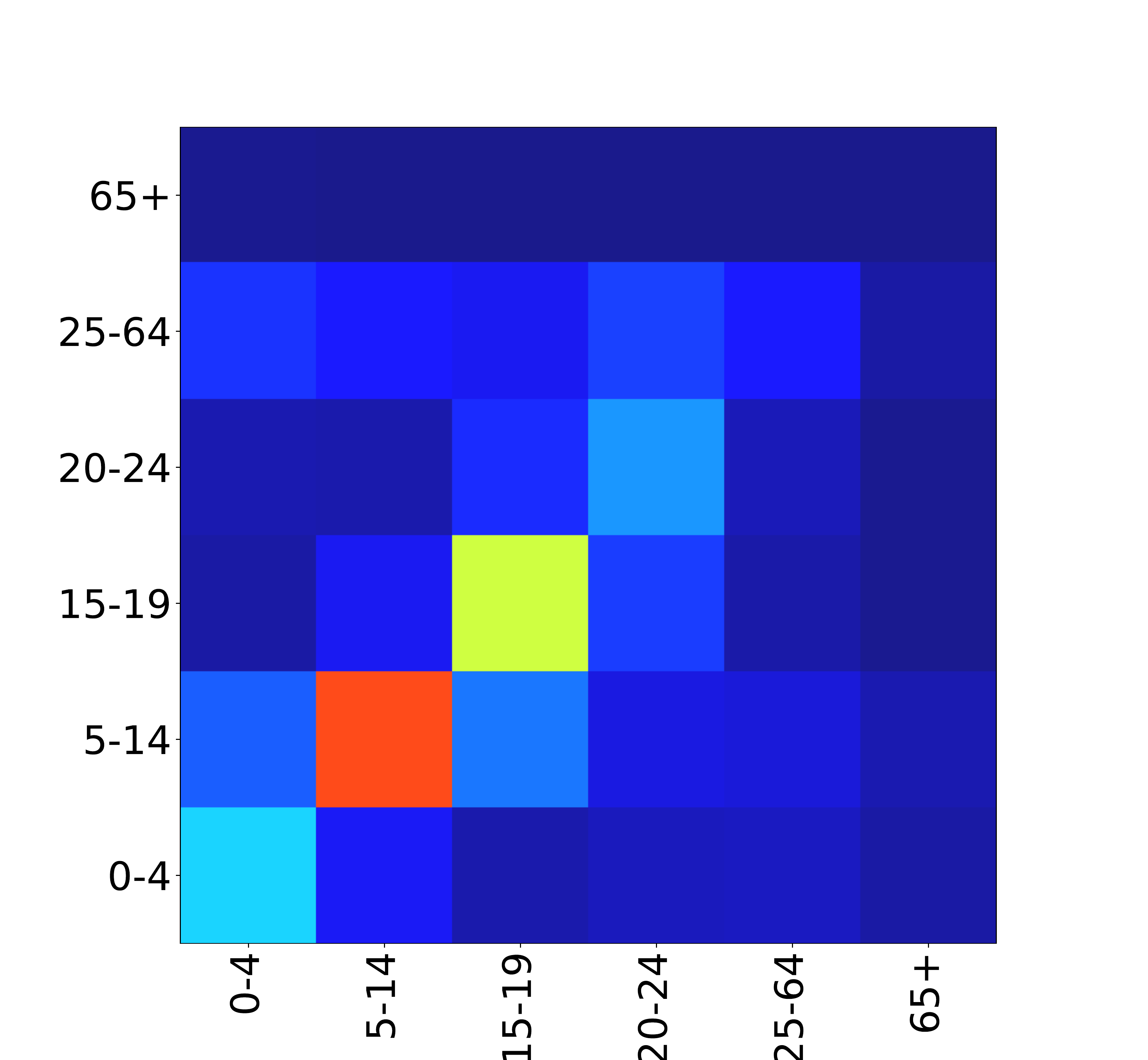}
         \caption{}
         \label{fig:other}
     \end{subfigure}
     \hfill
     \begin{subfigure}[b]{0.19\textwidth}
         \centering
         \includegraphics[width=\textwidth]{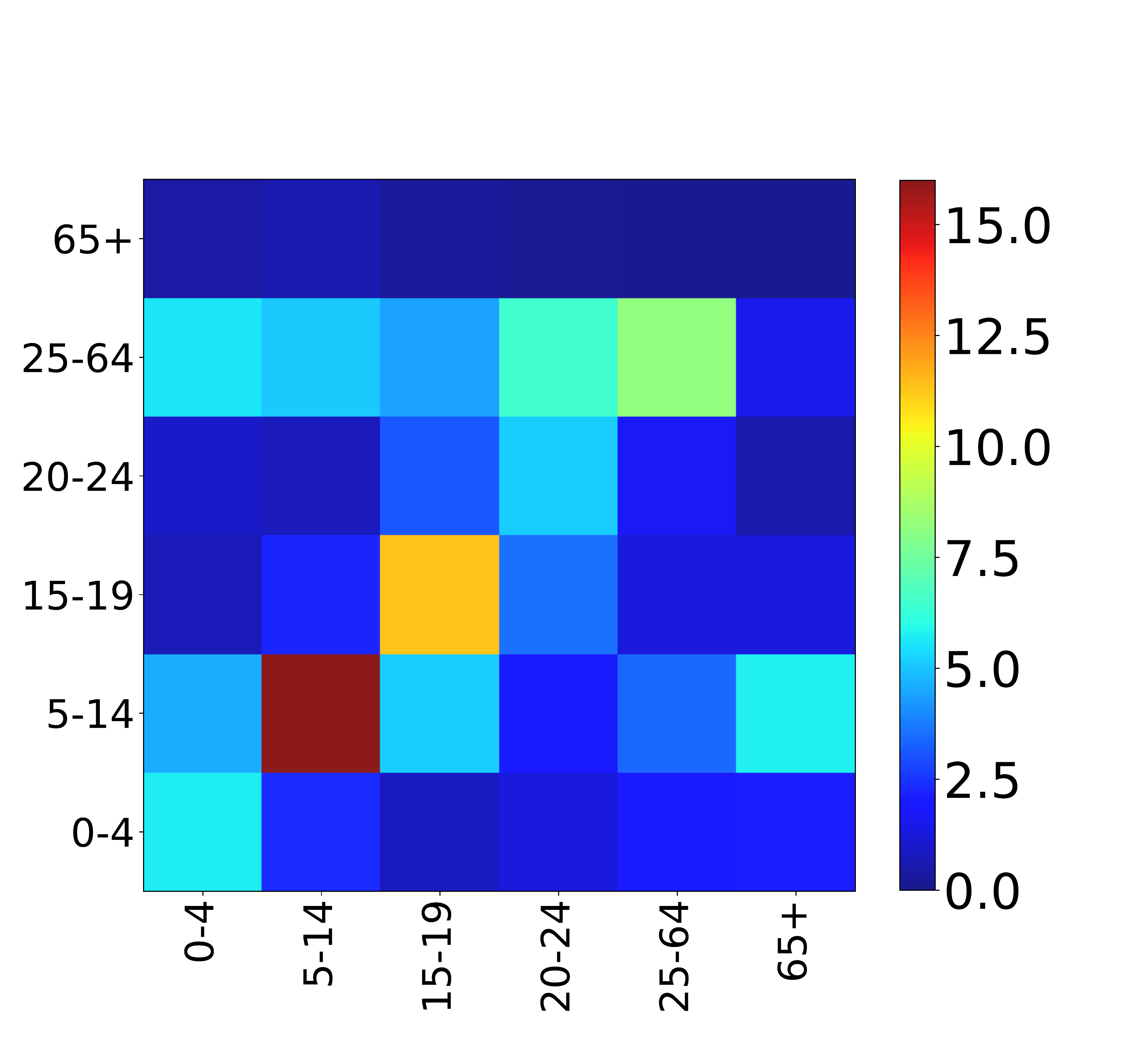}
         \caption{}
         \label{fig:full}
     \end{subfigure}  
        \caption{Heat maps of age-specific contact matrices in Kenya are shown. Estimated contact matrices at home (a), schools (b), workplaces (c), communities and other locations (d) and (e) is the total contact matrix obtained as the sum of the matrices in the four environments (home, school, at work and other locations). The ages of the participants are shown on the horizontal axis, while the ages of the reported contacts are shown on the vertical axis. The intensity of an entry color corresponds to the size of the contact rate between the interacting age groups.}
        \label{fig: contacts}
\end{figure}

Figure \ref{fig: contacts} shows contact matrices obtained for Kenya to demonstrate contact matrices in other African countries. 
Starting with the household setting in (Fig. \ref{fig: home}), we observe contacts between different age groups. There is strong contact between the elderly population and the young age group. There is also household contact between the working population (parents) and the young people (children). These contacts are generally lacking for the other age groups. Contacts at school (Fig. \ref{fig: school}) shows the highest frequencies of contact between children and young adults attending primary and secondary schools together. The interaction in tertiary education (20-24) and with instructors and other staff members in school is minimal.
The workplaces in Fig. \ref{fig:work} make it clear that most contacts occur between people of working age (25-64 years) and with people in tertiary education. This pattern can depend on local retirement, employment regulations, and culture in some countries.
According to Fig.\ref{fig:other}, the contacts in different places show strong connections along the diagonal for the young people (associative pattern), indicating a high population of young individuals in the region. 
 In Fig.\ref{fig:full}, contacts from all settings are reported; for interactions at home, interactions among students and communications from other places are represented on the main diagonal; The other references along the center are from the working-age population. The 25-64 age group has socialized more than other age groups due to their population size.

\subsection{Mathematical Model with heterogeneous mixing patterns}

\begin{figure}[ht]
    \centering
    \includegraphics[width=1\textwidth]{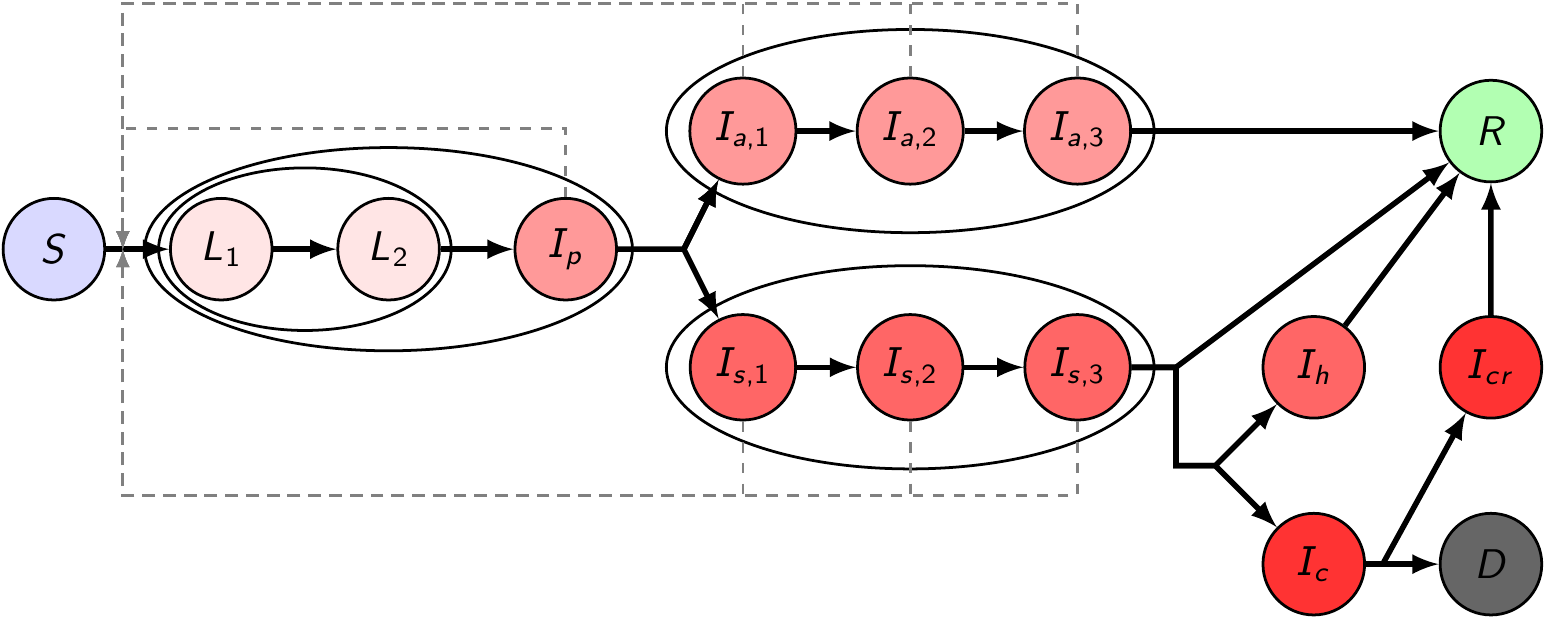}
    \caption{Flowchart of the transmission model from \cite{RB} used for demonstrating the proposed framework.}
    \label{img: diag}
\end{figure}

In this model (Fig. \ref{img: diag}), the population is divided into fifteen compartments: Susceptible ($S$) individuals who are at risk of infection, Latents $\left(L_{1}, L_{2}\right)$ who are already infected, but still are are not contagious. Presymptomatic $\left(I_{p}\right)$ people who are contagious but show no symptoms. Individuals from this compartment may progress to an asymptomatic with the following stages $\left(I_{a,1}, I_{a,2}, I_{a,3}\right)$ (mild symptoms) or symptomatic compartment with stages $\left(I_{s,1}, I_{s,2}, I_{s,3}\right)$ (severe symptoms). All those with mild symptoms recover and are collected in $R$. 
People in $I_s$ compartments can either recover or be hospitalized $\left(I_h\right)$. Some hospitalized individuals may require Intensive Care $\left(I_c\right)$. Individuals in the $I_c$ class can either switch to the $D$ class or proceed to $I_{\mathrm{cr}}$ before fully recovering. Each compartment has an age structure setup and is modeled with age-specific characteristics assumed to be the same for each country, except for transmission rates.

The age-specific model with contact matrices from selected African countries incorporated in \cite{PK} is given by the system of ordinary differential equations:

\begin{align}
 {S^i}'(t)={} & -\beta_0 \frac{S^i(t)}{W^{i}}\cdot\sum_{j=1}^{6} \mu_C^{(j,i)}\left[I_{p}^{j}(t) + \mathrm{inf}_a \sum_{k=1}^3 I_{a,k}^{j}(t) + \sum_{k=1}^3 I_{s,k}^{j}(t)\right] \nonumber\\ 
 {L_1^i}'(t)={} & \beta_0 \frac{S^i(t)}{W^{i}}\cdot\sum_{j=1}^{6} M^{(j,i)}\left[I_{p}^{j}(t) + \mathrm{inf}_a \sum_{k=1}^3 I_{a,k}^{j}(t) + \sum_{k=1}^3 I_{s,k}^{j}(t)\right] - 2 \phi_l L^i_1(t) \nonumber\\
 {L_2^i}'(t)={} & 2 \phi_l L_1^i(t) - 2\phi_l L_2^i(t),\nonumber\\
 {I_p^i}'(t)={} & 2 \phi_l L_{2}^{i}(t) - \phi_{p} I_{p}^{i} (t)\nonumber\\
 {I_{a,1}^i}'(t)={} & \theta^{i} \phi_{p} {I}_{p}^{i} (t) - 3 \delta_{a} I_{a,1}^{i}(t)\nonumber\\
 {I_{a,2}^i}'(t)={} & 3\delta_{a} I_{a,1}^{i}(t)- 3\delta_{a} I_{a,2}^{i}(t)\nonumber\\ 
 {I_{a,3}^i}'(t)={} & 3\delta_{a} I_{a,2}^{i}(t)- 3\delta_{a} I_{a,3}^{i}(t)\\ 
 {I_{s,1}^i}'(t)= {}& (1 - \theta^i) \lambda_{p} I_{p}^{i} - 3 \delta_{s} I_{s,1}^i(t)\nonumber\\
 {I_{s,2}^i}'(t)= {} & 3 \delta_{s} I_{s,1}^i(t) - 3 \delta_{s} I_{s,2}^i(t)\nonumber\\
 {I_{s,3}^i}'(t)= {} & 3 \delta_{s} I_{s,2}^i(t) - 3 \delta_{s} I_{s,3}^i(t)\nonumber\\
 {I_h^i}'(t)={} & h^i (1 - \zeta^i) 3 \delta_{s} I_{s,3}^i(t) - \delta_h I_h^i(t)\nonumber\\
 {I_c^i}'(t)={} & h^i \zeta^i 3 \delta_{s} I_{s,3}^i(t)-\delta_c I_c^i(t)\nonumber\\
 {I_{\mathrm{cr}}^i}'(t)={} &  (1 - \eta^{i}) \delta_{c} I_{c}^{i} (t) -\delta_{\mathrm{cr}} I_{\mathrm{cr}}^{i} (t)\nonumber\\
 {R^i}'(t)={} & 3 \delta_{a} I_{a,3}^{i} (t) + (1- h^{i}) 3 \delta_{s} I_{s,3}^{i} (t) + \delta_{h} I_{h}^{i} (t) + \delta_{\mathrm{cr}} I_{\mathrm{cr}}^{i} (t)\nonumber\\   
 {D^i}'(t)={} & \eta^i \delta_c I_c^i(t),\nonumber
 \label{eq:model}
\end{align}
where the age groups are represented by the index $i \in 1,..., 6$, and $\beta_0$ is the likelihood of transmission given a contact. The upper index $i$ for each of the listed parameters denotes the relevant age group. One in every $\left(1-\theta^{i}\right)$ exposed individuals will experience symptoms throughout their illness, compared to $\left(\theta^{i}\right)$ who will not. The incubation period lasts an average of $\left(\phi_{L,1}^i\right) ^{-1} + \left(\phi_ {L,2}^i\right)^{-1} + \left(\phi_{I,p}^i\right)^{-1}$ days. Similarly, infected individuals in asymptomatic and symptomatic compartments will spend an average of $\left(\delta_{a,1}^i\right)^{-1} + \left(\delta_{a,2}^i\right)^{-1} + \left(\delta_{a,3}^i\right)^{-1}$ and $\left(\delta_{s,1}^i\right)^{-1} + \left(\delta_{s,2}^i\right)^{-1} + \left(\delta_{s,3}^i\right)^{-1}$ days respectively. $h^{i}$ of the infectious compartment $\left(I_{s,3}^{i}\right)$ will require hospitalization, while the other $1-h^{i}$ will recover on their own. A portion of individuals, $\zeta^{i}$ who require hospitalization require intensive care. The average length of stay for the hospitalized classes $I_h^{i}, I_c^{i}$, and $I_{cr}^{i}$ is represented as $\left(\delta_{h}^i\right)^{-1}$, $\left(\delta_{c}^i \right)^{-1}$, and $\left(\delta_{cr}^i\right)^{-1}$ respectively. A portion $\left(1-\eta^{i}\right)$ will move on to the $I_{cr}^{i}$ class while $\eta^{i}$ will pass away from the illness. The relative infectiousness of $I_a$ in comparison to $I_s$, is shown by the symbol $\mathrm{inf}_a$. For more details about the other parameters and the methodology for parametrization, see \cite{Bra, Kim, RB, ZE} and Table \ref{table:params} above.
Fig. \ref{img: diag} shows the transmission dynamics of our model for one age group.

The basic reproduction number, $\mathcal{R}_{0}$, is calculated using the 
Next Generation matrix method $(\mathrm{NGM})$. It is the dominant eigenvalue of the $\mathrm{NGM}$ method \cite{DH, EZ, KK, KR, MX, RB},
$$\mathrm{NGM} = -\beta_0\cdot T \cdot\Sigma^{-1}$$ 

It can easily be verified using the approach that the $\mathcal{R}_{0}$ for system of ordinary differential equations satisfies

$$\mathcal{R}_0 = \beta_0 \cdot \rho \left(-T  \cdot\Sigma^{-1}\right).$$

For a constant value of $\mathcal{R}_{0}$ we can determine the transmission rates of the countries. \cite{IR} Estimated basic reproduction number of African countries using the $\mathrm{SEIR}$ model and a Bayesian inference approach. The authors found an average and median $\mathcal{R}_{0}$ of 3.68 and 3.67, respectively. Thus, a $\mathcal{R}_{0}$ of 3.68 is used to generate base transmission rates $\beta_{0,C}$ for each country $C\in\{1, 2, \dots, 32\}$. With these calculated country-specific transmission rates, we can scale all country contact matrices  to make them comparable,
\begin{equation*}
    S_{C} = \beta_{0,C} \cdot \mu_{C}.
\end{equation*}

Countries can then be compared using these standardized contact matrices $S_{C}$.

\subsection{Dimensionality reduction}

Modern applications of statistical theory involve vast amounts of data with a more significant number of features compared to data points which pose challenges to optimization algorithms and training performance in terms of computational speed and computer memory. To address these issues, data scientists apply feature selection or dimensionality reduction techniques \cite{AK, HF}. One of the techniques, classical principal components analysis (PCA), finds an optimal projection matrix that maps each point onto a low-dimensional space. The covariance method and singular value decomposition can effectively determine the principal components. PCA was applied to the socioeconomic data to reduce the dimensionality of the 28 features for each country into four feature vectors;

\begin{equation}
    \tilde{F}\in\mathbb{R}^{32\times 4}.
    \label{eq: pca}
\end{equation}

\begin{figure}{}
     \centering
     \begin{subfigure}[b]{0.24\textwidth}
         \centering
         \includegraphics[width=\textwidth]{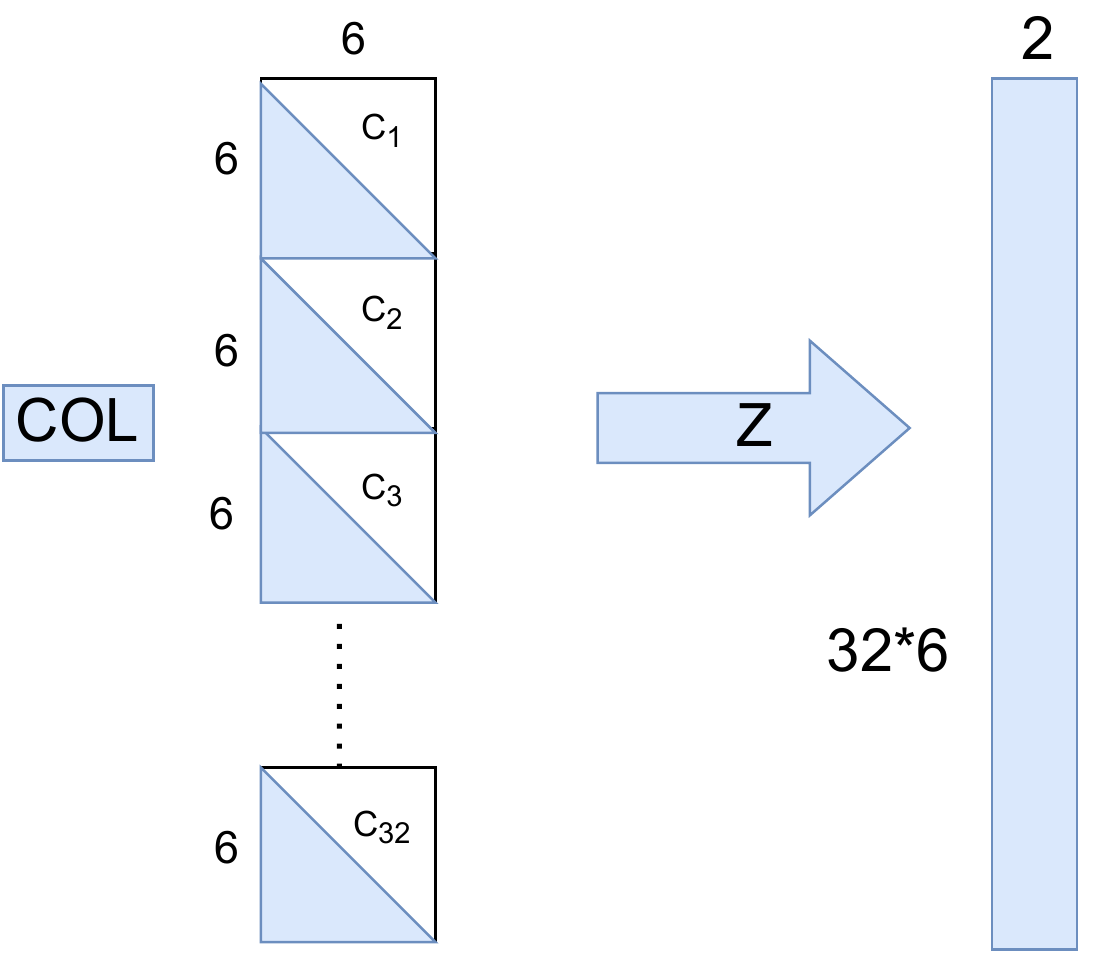}  
         \caption{}
         \label{fig: col}
     \end{subfigure}
     \hfill
     \hspace{0.6cm}
     \begin{subfigure}[b]{0.24\textwidth}
         \centering
         \includegraphics[width=\textwidth]{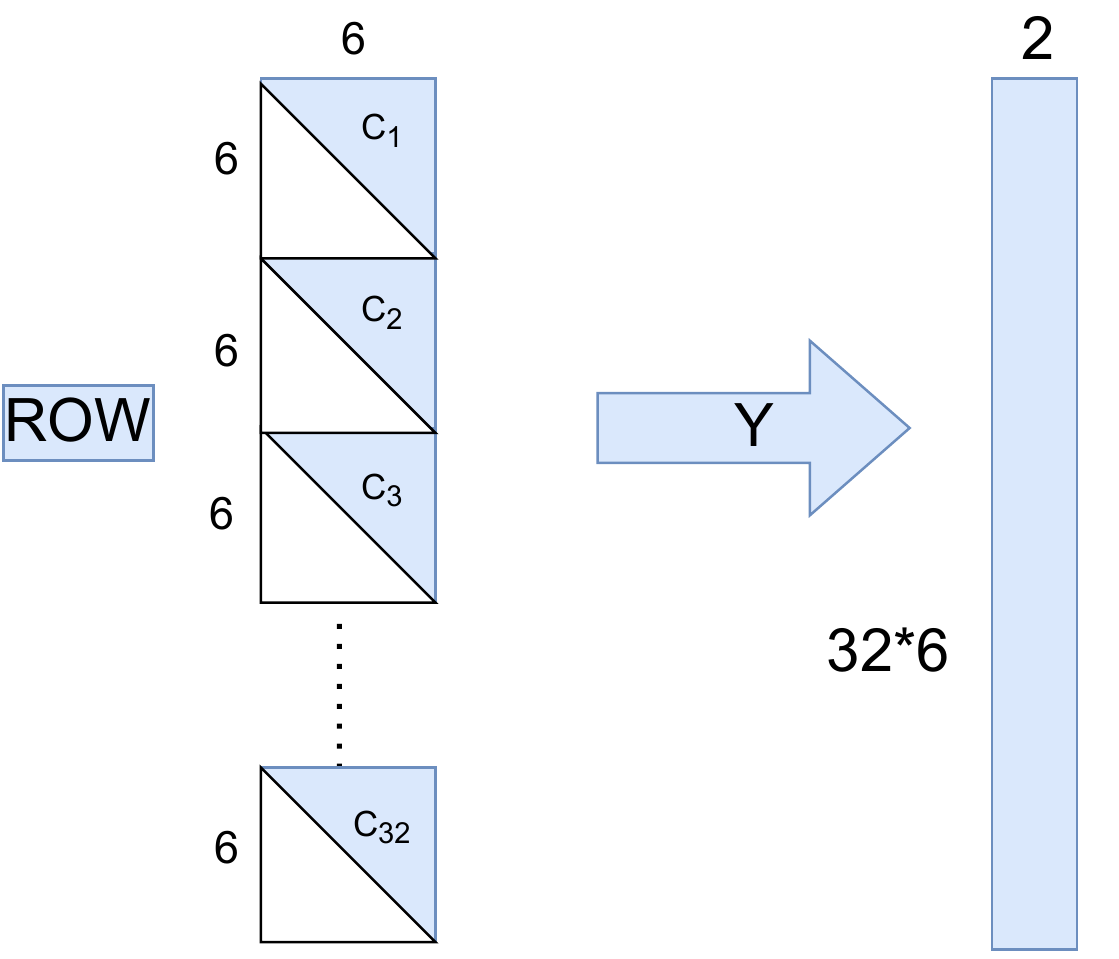}
         \caption{}
         \label{fig: row}
     \end{subfigure}
      \hspace{0.6cm}
      \hfill
     \begin{subfigure}[b]{0.24\textheight}
         \centering
         \includegraphics[width=\textwidth]{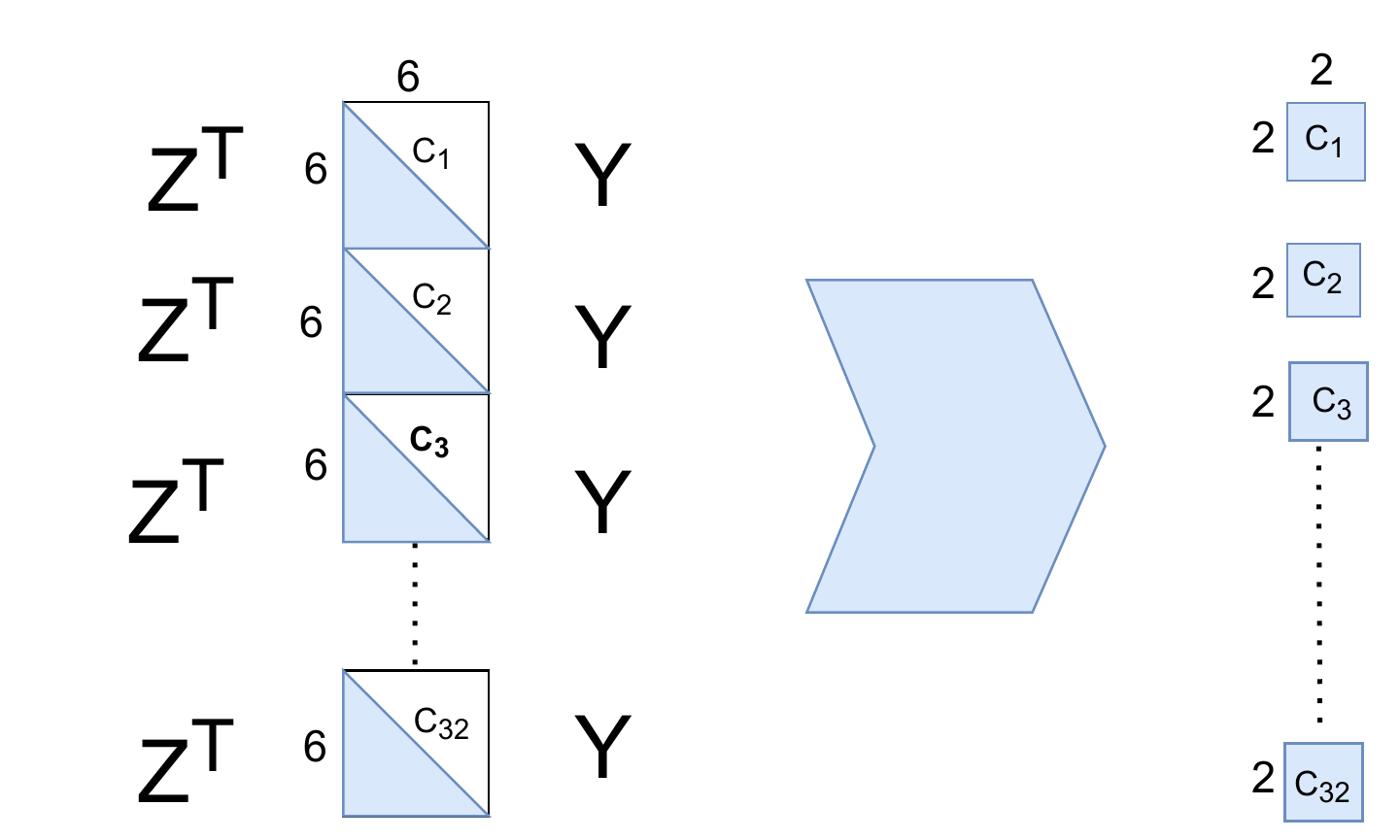}
         \caption{}
         \label{fig: both}
     \end{subfigure} 
        \caption{The steps for dimensionality reduction utilizing the $(2D)^ 2 PCA $ are displayed. The concatenated matrix is initially projected to 2D while taking into account the column direction (a). By projecting the concatenated transposed matrices, two principal components were also preserved in the row direction (b). In (c), we project each of the $S_C$ matrices using the two projection matrices $Y$ and $Z$ to obtain the reduced matrix of dimension $2 \times 2$.}
        \label{fig: dim}
\end{figure}

The feature vectors from the contact matrix will be incorporated into these feature vectors from the PCA.
Since we have high-dimensional matrix-valued data, $M \in \mathbb{R}^{32 \times (6 \times 6)}$, dimensionality reduction techniques will be fundamental to improving learning performance. When using classical PCA on matrix-valued data, it is challenging to precisely calculate the covariance matrix and the eigenvectors since it concatenates 2D matrices into 1D vectors, creating a high-dimensional vector space. Thus, it is crucial to use a 2-Directional 2-Dimensional PCA ($(2D)^2 PCA $) reduction technique that is precise and scales well for high-dimensional data \cite{EZ, WS, ZT}. The approach is applied to country-specific complete contact matrices.
Using this approach, we first concatenate each country-specific complete contact matrices row by row, as shown in Fig. \ref{fig: dim}, to obtain the
matrix $X_{\mathrm{col}}=[S_1^\top, S_2^\top, \dots, S_{32}^\top]^\top \in \mathbb{R}^{(32 \cdot 6) \times 6}$. In the column direction, we projected $X$ onto the plane spanned by the columns of $Z \in \mathbb{R}^{6 \times 2} $, resulting in $$\hat{X}_{\mathrm{COL}} = X_{\mathrm{COL}} \cdot Z \in \mathbb{R}^{(32 \cdot 6) \times 2}.$$

Considering the row direction as in Fig. \ref{fig: row}, we projected concatenation of the transposed matrices $X_{\mathrm{ROW}} = [S_{1}, S_{2}, \dots, S_{32}]^\top \in \mathbb{R}^{(32 \cdot 6) \times 6}$ onto the plane spanned by the columns of $Y \in \mathbb{R}^{6 \times 2} $ to obtain, $$\hat{X}_{\mathrm{ROW}} = X_{\mathrm{ROW}} \cdot Y  \in \mathbb{R}^{(32 \cdot 6) \times 2}.$$

Given the two projection matrix $Y$ and $Z$ in Figure \ref{fig: both}, we project the country-specific matrices $S_C$ to obtain the smaller matrices,
$$\hat{S}_C = Z^\top \cdot S_C \cdot Y \in \mathbb{R}^{2 \times 2},$$ where $C\in\{1, 2, \dots, 32\}$.
The smaller matrices were then transformed to derive the feature vectors,

\begin{equation}
    s_{C}\in\mathbb{R}^4.
    \label{eq: dpca}
\end{equation}

The feature vectors in Eq. \ref{eq: dpca} and Eq.\ref{eq: pca} are appended column by column to produce the appended feature vectors $m_{C}\in \mathbb{R}^8$. 
The countries are clustered using this appended feature vectors, $m_{C}$.

\subsubsection{Clustering analysis}

Clustering has become standard due to the large amounts of data available in different companies. Clustering involves grouping data into clusters such that the similarity within the set is high and between the collections is low. It can be carried out in various methods, such as feature selection, density-based clustering, probabilistic clustering, grid-based clustering, and spectral clustering, among other techniques \cite{AK, BB, HF}. Distance-based techniques are commonly used because they can work with any data based only on an appropriate distance measure. We used the hierarchical clustering technique because we have small data sets, and the method is widely used in the literature. This technique was used in \cite{EZ, Mis} to group countries and regions based on contact matrices.

The technique creates clusters by developing a tree-like structure called a dendrogram. The approach can be achieved by considering bottom-up (agglomerative) or top-down (divisive) modes. We considered the agglomerative method because it is easy to implement the algorithm as it builds the clustering by treating the data points as single clusters in the first step. The sets are then merged at each step until the final merge containing all data points is obtained. The groups of clusters are connected using a linkage method; single, average, ward, or complete \cite{CC, ML, NB, RU, SC}. 

The study considered the complete method as it preserves the structure of the clusters, produces compact groups, and is less sensitive to noise and outliers \cite{AK, ML}.
The method considers the longest inter-cluster distance among the observations as the inter-cluster distance. In each iteration, it tries to minimize this distance. For two clusters $C_{i}$ and $C_{j}$ the method can be defined as follows:

\begin{equation*}
     D_{\mathrm{CL}}(C_{i}, C{_j}) = \sum_{x \in C_{i}, y \in C_{j} } d_{dist}(x, y).   
\end{equation*}
where $dist$ is the chosen distance measure. This study considers Euclidean distance, but other distance metrics can be used to reach desired clusters. After creating the dendrogram, it is easy to determine the number of groups. As a rule of thumb, we consider the vertical distance with the most significant distance between the two closest merges to achieve desired sets without rerunning the clustering procedure.

\begin{figure}[ht]
    \centering
    \includegraphics[width=1\textwidth]{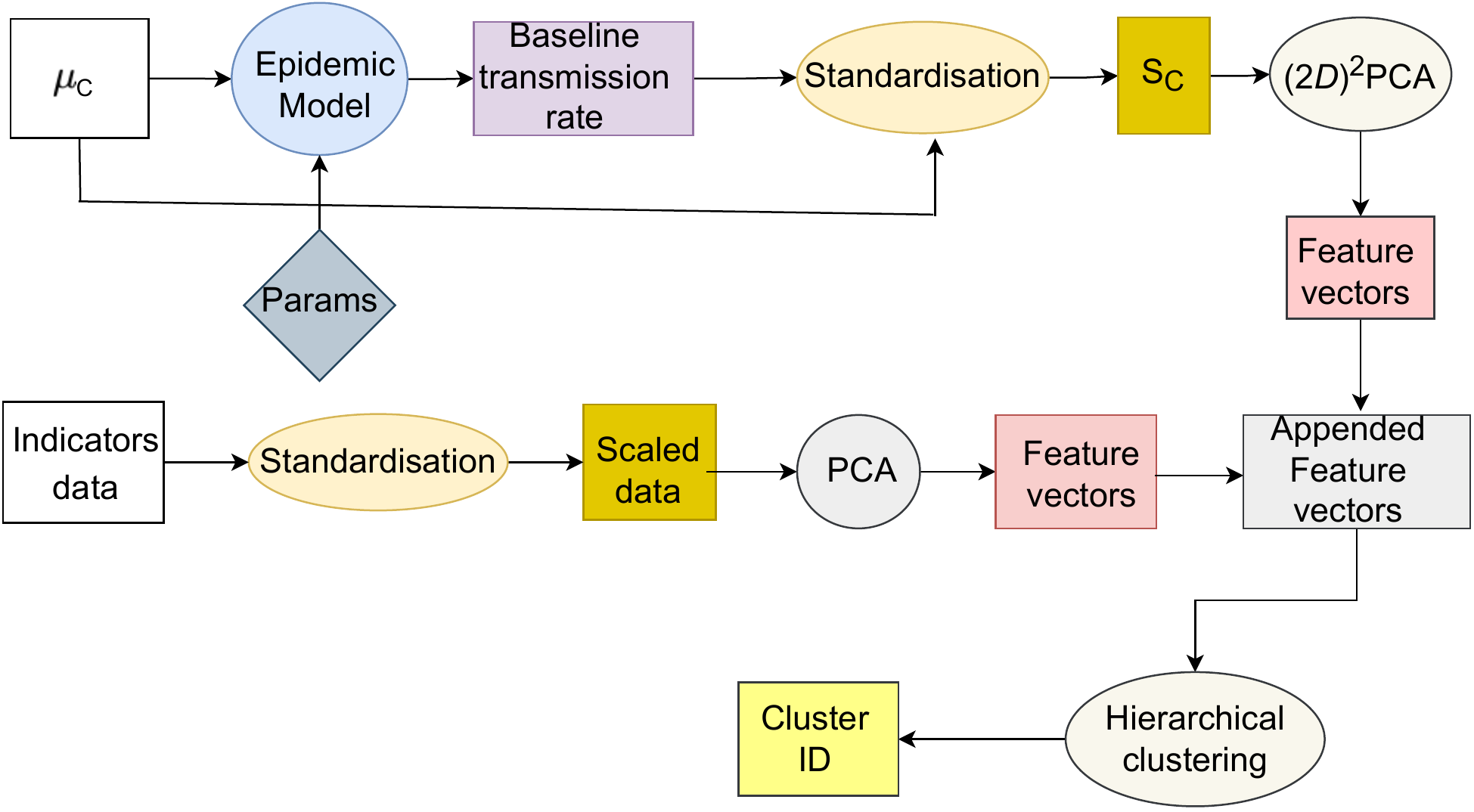}
    \caption{Flowchart showing a summary of the research design. Rectangles show the results of a method. As can be seen, the contact matrices and socioeconomic indicators are the input of the study. We apply standardization to both inputs to get the scaled data. In the next steps we apply PCA to the scaled data of indicators and $(2D)^2 PCA$ to scaled contact data. The output feature vectors are then merged to get an appended one. Then hierarchical clustering is performed to get the clusters.}   
    \label{fig: method}
\end{figure}

\section{Results}

We analyzed synthetic contact matrices from \cite{PK} for 32 countries in Africa and integrated socioeconomic data of the countries available in \cite{Z}.

\subsection{Visualization of the countries: socioeconomic similarities}

\begin{figure}{}
     \centering
     \begin{subfigure}[b]{0.49\textwidth}
         \centering
         \includegraphics[width=\textwidth]{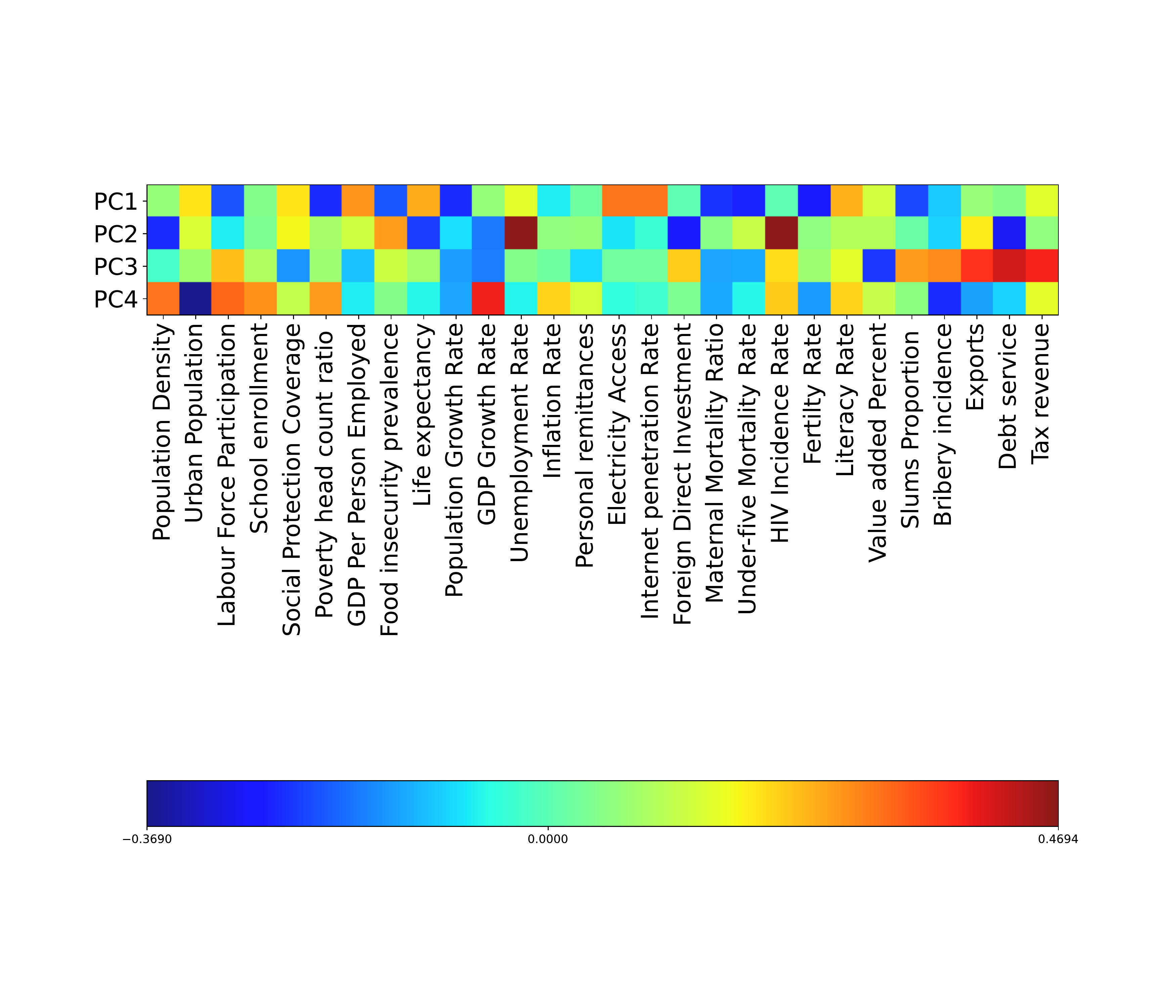}
         \caption{}
         \label{fig: components}
     \end{subfigure}
     \hfill     
     \begin{subfigure}[b]{0.49\textwidth}
         \centering
         \includegraphics[width=\textwidth]{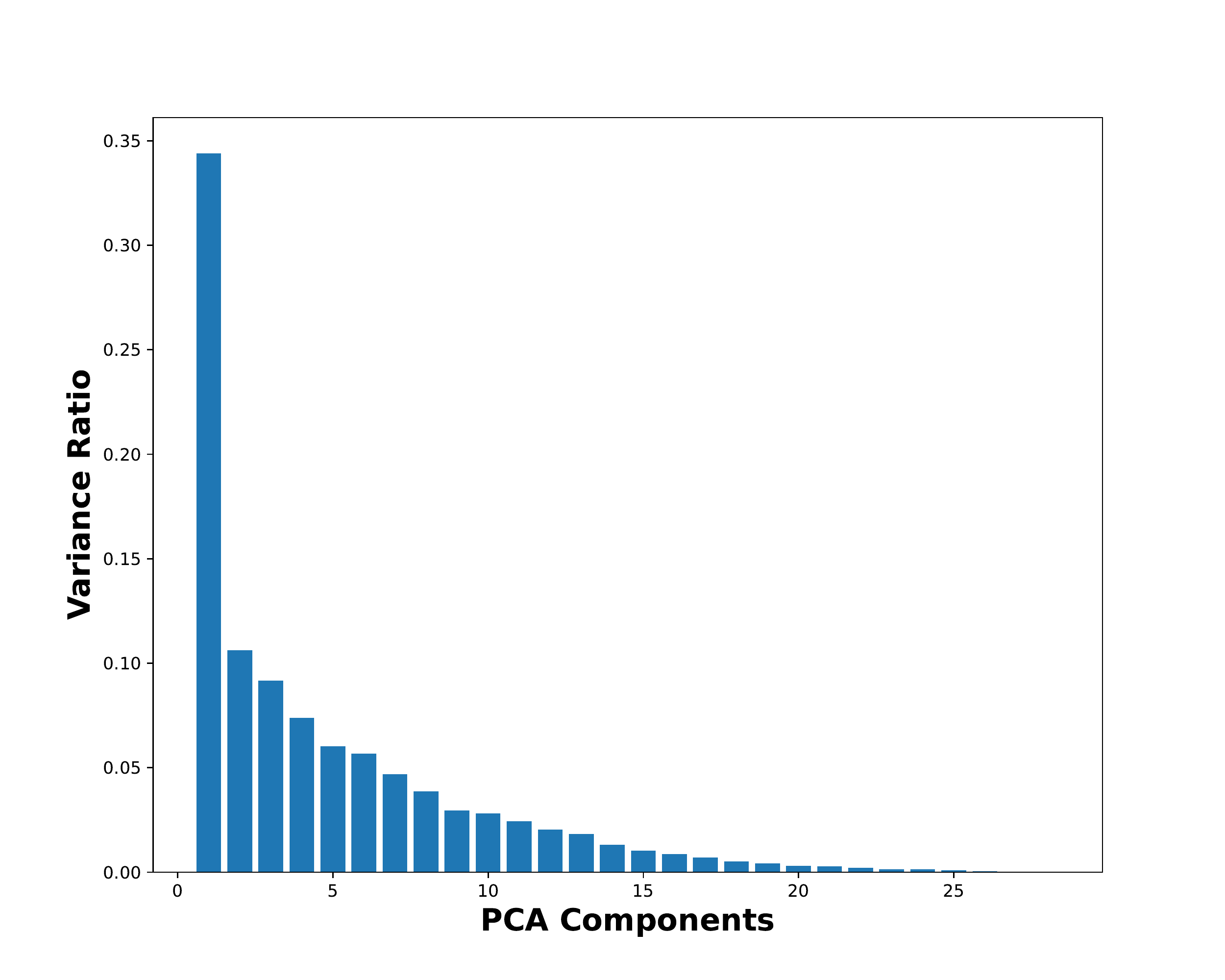}
         \caption{}
         \label{fig: var}
          \end{subfigure} 
          \hfill
          \begin{subfigure}[b]{0.5\textwidth}
         \centering
         \includegraphics[width=\textwidth]{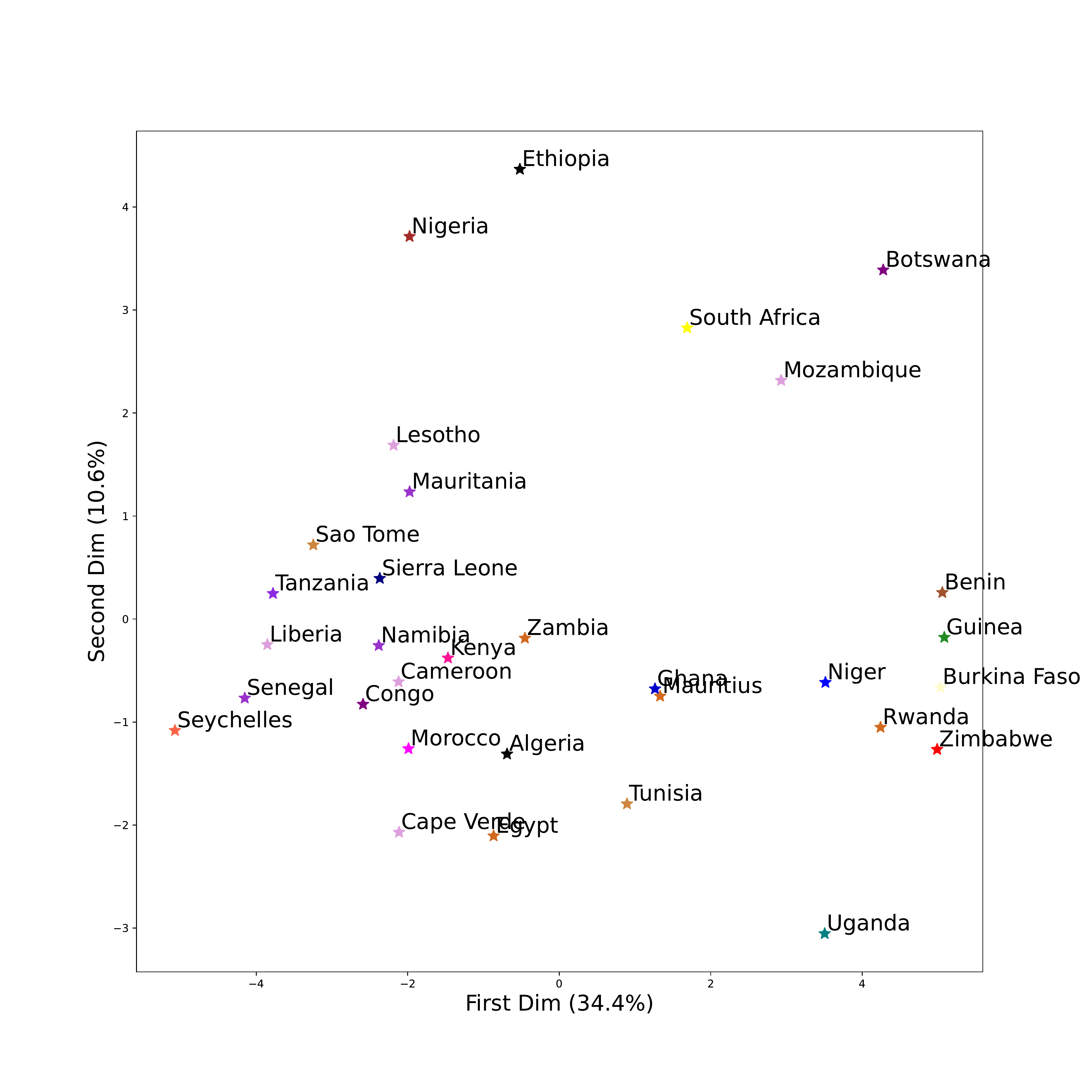}
         \caption{}
         \label{fig: countries}
     \end{subfigure}
        \caption{Part (a) presents the factor loadings for the first four principal components to indicate the correlation coefficient between the original feature and the principal component. Part (b) is a screen plot that plots eigenvalues or singular values and shows the proportion of the total variance represented by the principal components.  Part (c) displays the projection of countries onto a 2D space based on their socioeconomic factors. The two dimensions (first dim and second dim) explained 45 percent of the socioeconomic variance in the data given the 28 attributes used in this analysis. Countries are colored with different random colors. }
        \label{fig: proj}
\end{figure}

Using principal component analysis (PCA), we found that the first principal component (PCA) explained 34.4 percent (Fig. \ref{fig: var}), and the four principal components explained 61.6 percent of the total socioeconomic variance among the 28 variables considered (Figure \ref{fig: var}). From Figure \ref{fig: components}, the original socioeconomic variables that contributed most to the first principal components are internet penetration, access to electricity, and GDP per worker. The second component consists mainly of unemployment rates and HIV incidence rates. The third PC is dominated by exports, debt service, and tax revenue, while the fourth PC is the GDP growth rate component. The projection of the 32 African countries in two dimensions resulted in a roughly distributed scatter diagram in Fig. \ref{fig: countries}. From the projection, countries such as South Africa, Mozambique, and Botswana lie beside each other. The same applies to the North African countries (Algeria, Morocco, Tunisia, and Egypt). Zimbabwe, Burkina Faso, Guinea, and Benin show similarities in socioeconomic indicators. Uganda has been pushed to the lower side of the plot, indicating differences in socioeconomic indicators.

Using four principal components, the hierarchical clustering algorithm generated a dendrogram in Figure \ref{fig: socioecon} to represent clusters of countries based on the socioeconomic indicators. As a rule of thumb, we cut the dendrogram at the height of 8, which divides the countries into four groups (Fig. \ref{fig: socioecon}).

\subsection{Country clustering and analysis of social contact patterns and associated socioeconomic indicators.}

\begin{figure}{}
     \centering
     \begin{subfigure}[b]{0.49\textwidth}
         \centering
         \includegraphics[width=\textwidth]{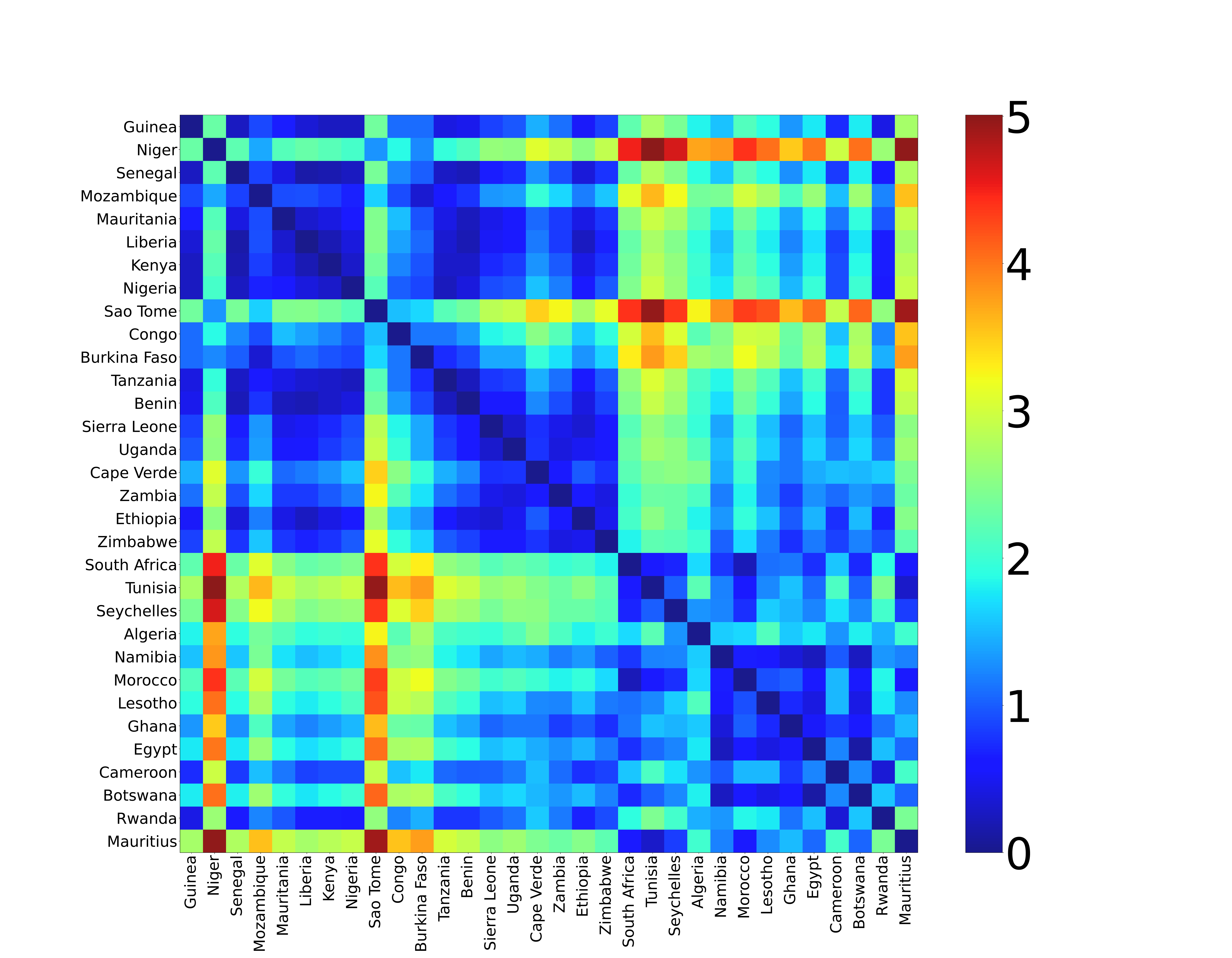}
         \caption{}
         \label{fig: dpca_contacts}
     \end{subfigure}
     \hfill
     \begin{subfigure}[b]{0.49\textwidth}
         \centering
         \includegraphics[width=\textwidth]{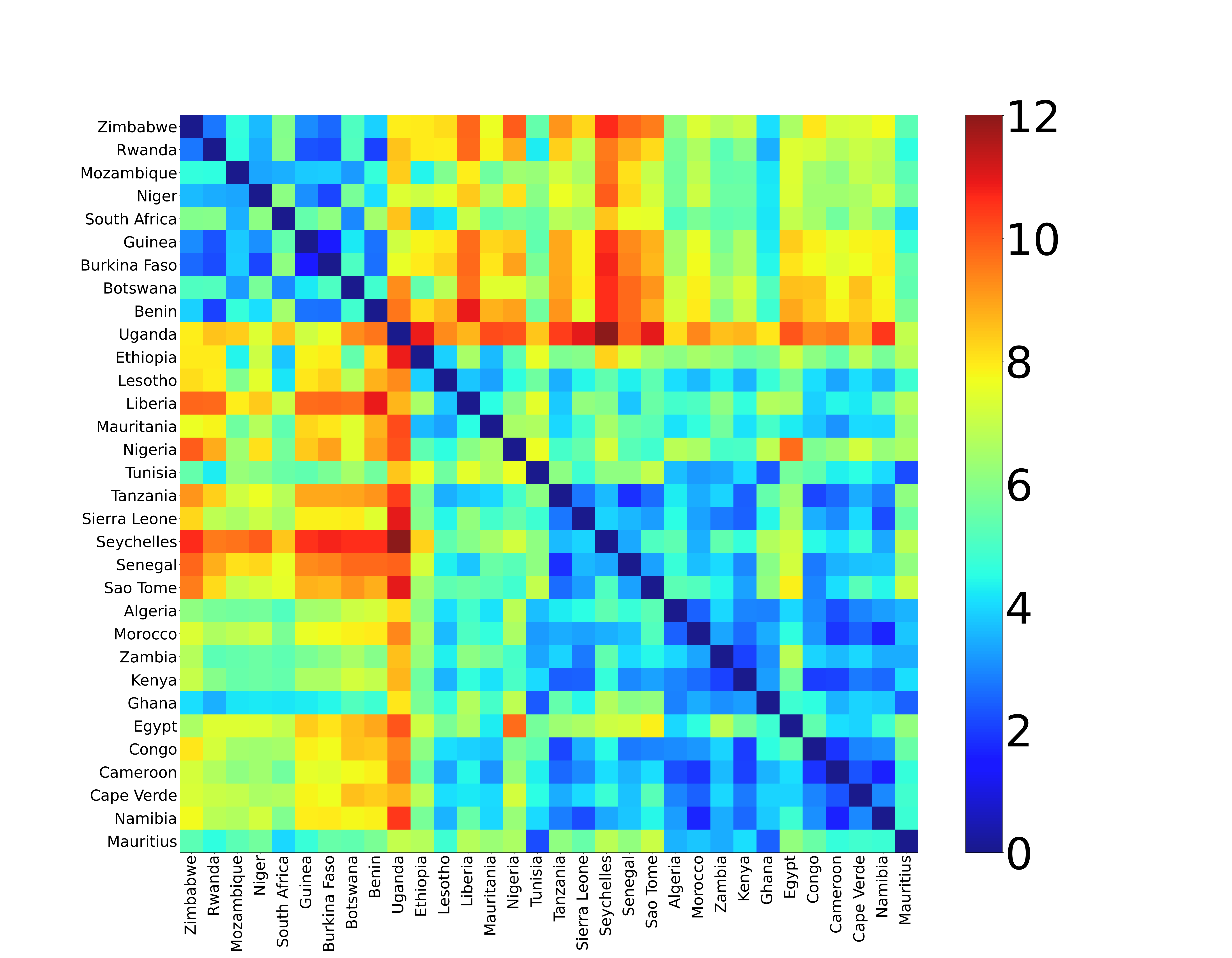}
         \caption{}
         \label{fig: dpca_combine}
     \end{subfigure}
        \caption{
        The Euclidean pairwise distance matrices calculated for all countries based on the feature vector from the $(2D)^2 PCA$ approach. Figure (a) corresponds to the pairwise distance based on the social contact matrices, while figure (b) gives the distances based on the contact matrices and socioeconomic indicators. In order to observe a pattern and present the elements in a nicer way, the dendrogram has been used to sort the rows and columns.}
        \label{fig: dpca_ord}
\end{figure}

By scaling the contact matrices and performing dimensionality reduction, the calculated feature vectors in Eq. \ref{eq: dpca} can then be integrated with the feature vectors from the socioeconomic indicators in Eq. \ref{eq: pca} to calculate the pairwise Euclidean distance for the countries. The dissimilarity matrix can then be permuted using the dendrogram to see a nicer plot and cluster structure. Figure \ref{fig: dpca_contacts} shows the distance matrix calculated from only the social contacts: as you can see, countries with dark blue color entries, such as Mauritania, Liberia, Kenya, Guinea, Ethiopia, among others, have a short paired distance, and such countries have closely related social contact patterns and are likely to belong to the same cluster if you consider only the social mixing patterns. The same applies to other clustering structures. Countries like Sao Tome and Niger have orange cell colors, with most countries showing their different way of socializing compared to other countries. Such countries may develop different intervention strategies.

Figure \ref{fig: dpca_combine} represents the distance matrix based on the social contacts and considering the country's socioeconomic indicators. We can see that most countries were influenced by socioeconomic performance. Uganda is more affected by socioeconomic variables. Similarly, Seychelles, Senegal, Liberia, Nigeria, and Tanzania appear to be affected, and this effect can be seen directly from the distance matrix. The permutation of the matrix provides a better presentation of the matrix.

\begin{sidewaysfigure}
  \centering
    \includegraphics[width=\textheight]{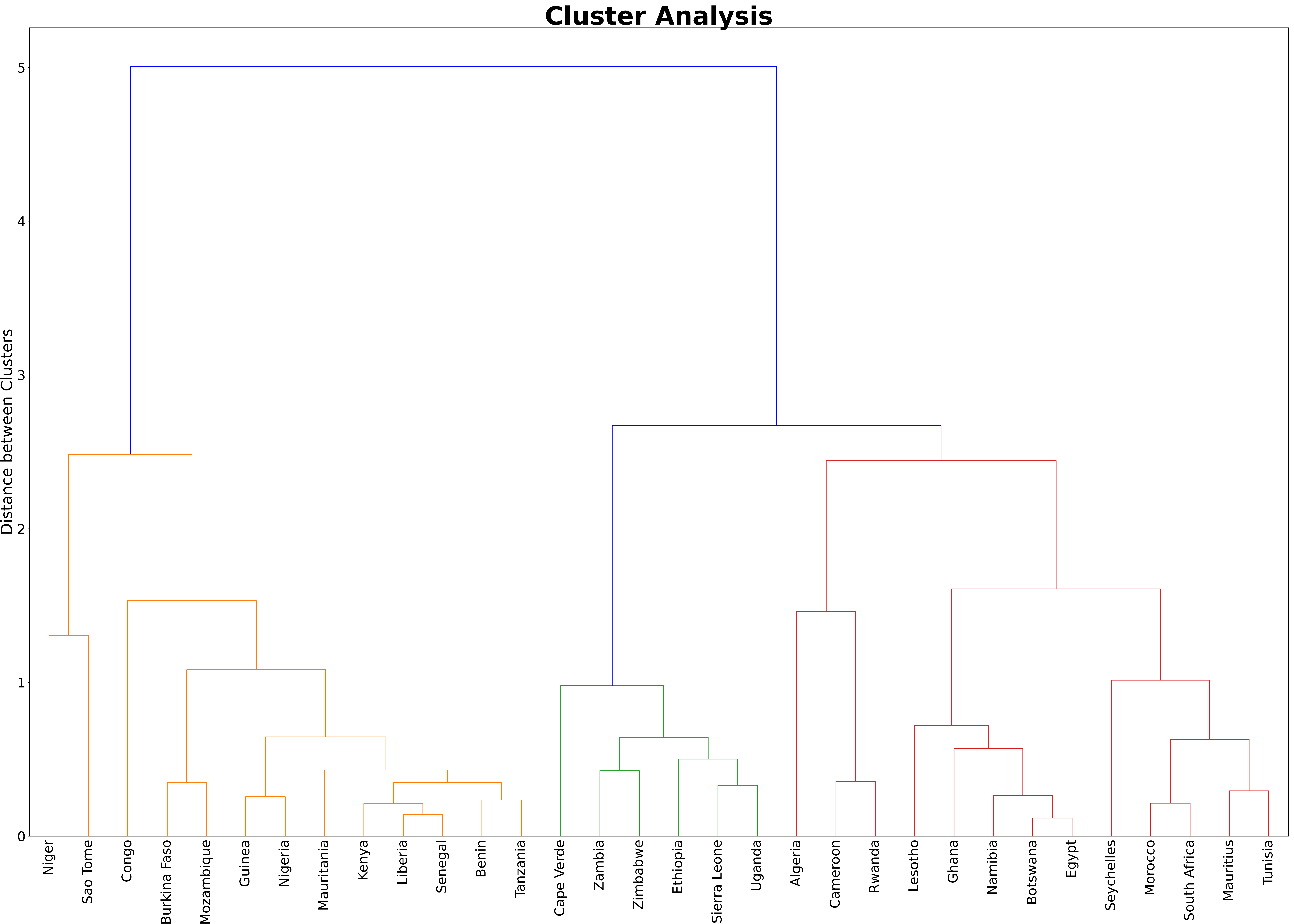}
    \caption{Dendrogram related to clustering of countries based on feature vectors of the contact matrices. The clustering is implemented with $(\mathrm{2D})^{2}$ PCA. Along the vertical axis we have the distance between the two clusters that are being merged. The clusters are agglomerated from bottom to top. Chopping down the tree at a height of 2.5 yields three clusters highlighted in orange, green, and red. Countries with the same hue belong to the same cluster.}
  \label{fig: dpca dendo conta}
\end{sidewaysfigure}

\begin{sidewaysfigure}
    \centering
    \includegraphics[width=\textheight]{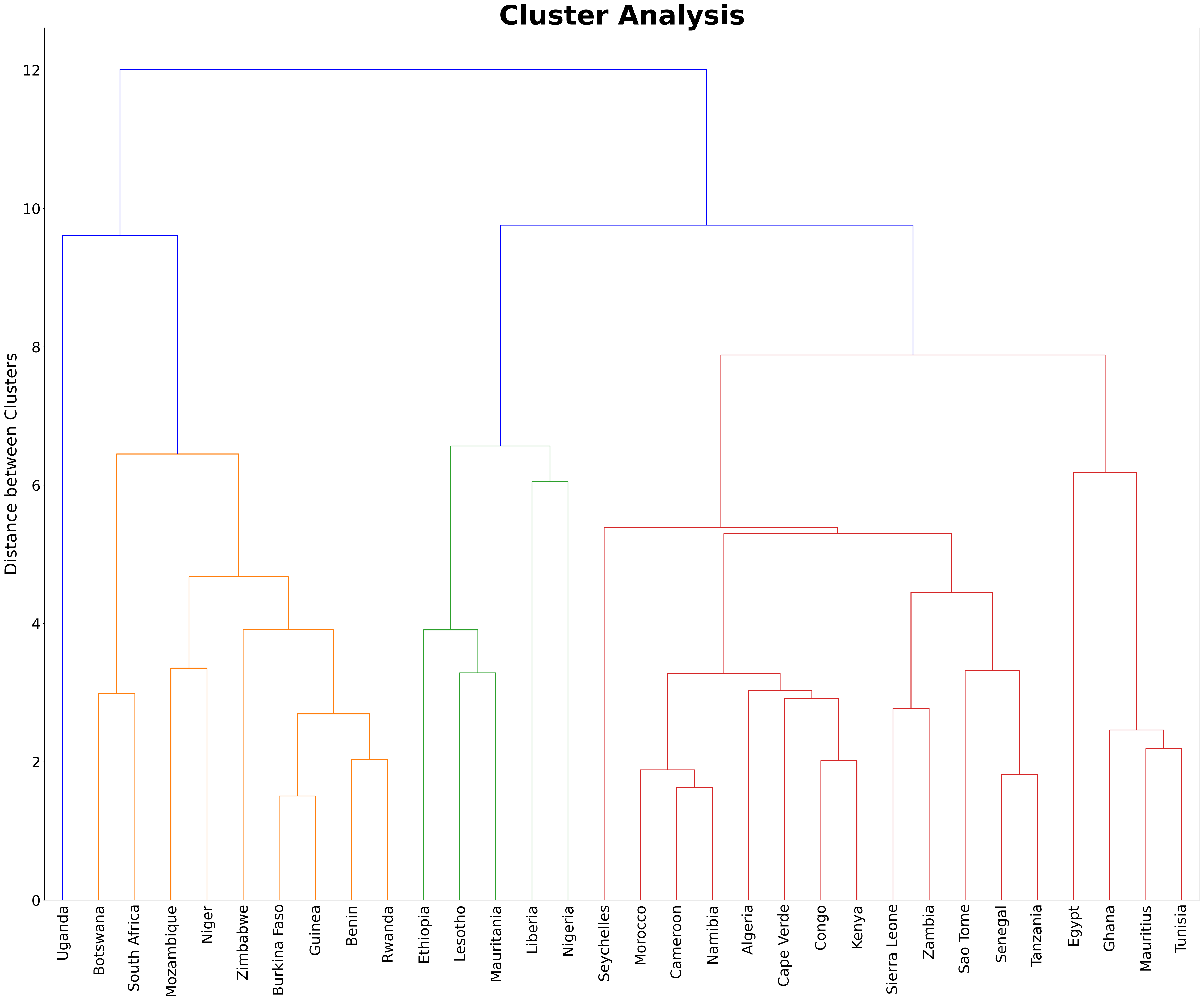}
    \caption{Dendrogram related to the grouping of countries based on feature vectors of contact matrices and socioeconomic indicators. The clustering is implemented with $(\mathrm{2D})^{2}$ PCA and PCA. Along the vertical axis we have the distance between the two clusters that are being merged. The clusters are agglomerated from bottom to top. Chopping down the tree at a height of 8 yields four clusters highlighted in orange, green, blue and red. Countries with the same hue belong to the same cluster.}
  \label{fig: dpca cont socio}
  \end{sidewaysfigure}

\begin{sidewaysfigure}
\centering
\captionsetup{width=0.8\textwidth} \includegraphics[width=\textheight]{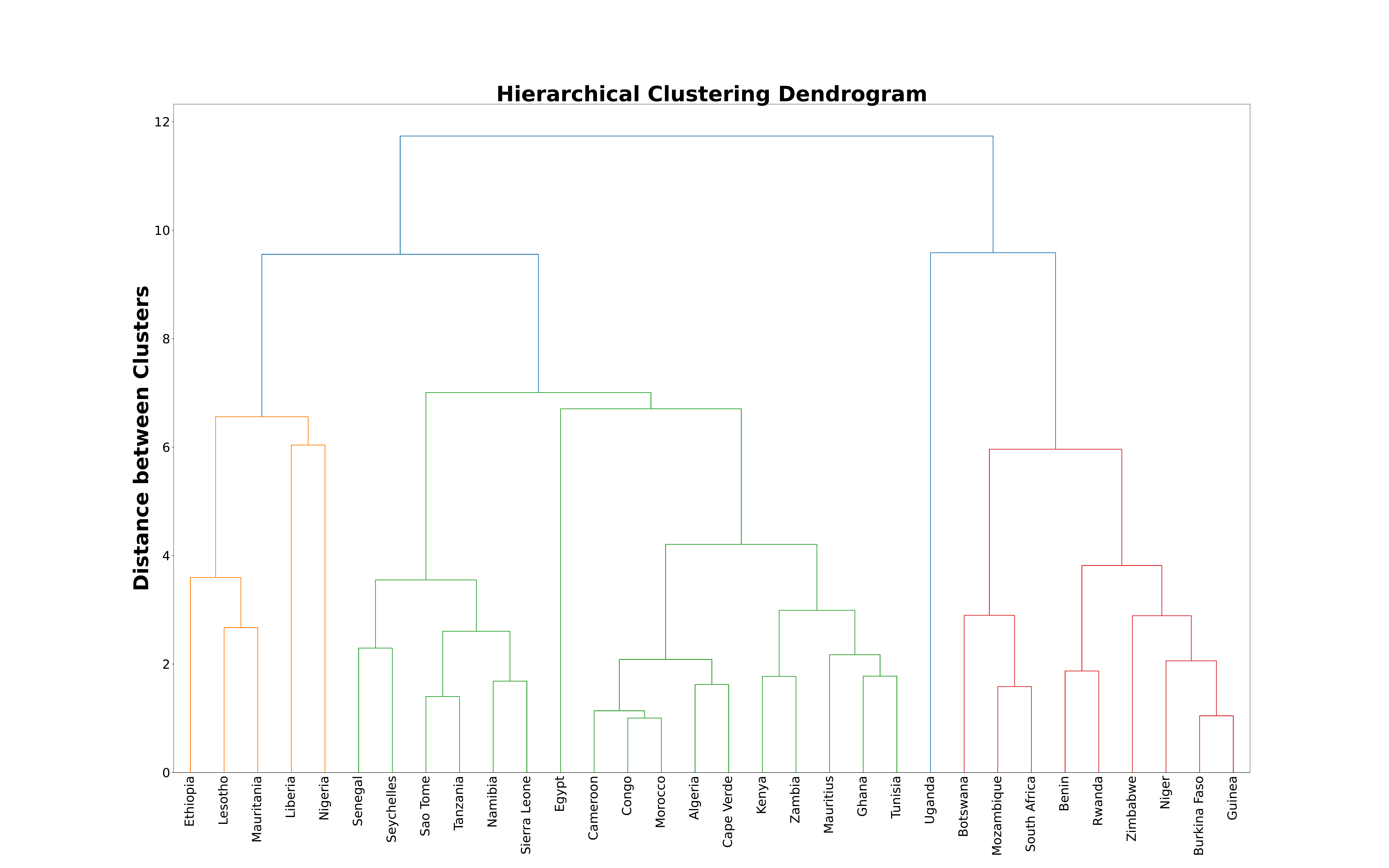}
    \caption{Dendrogram associated with the clustering of countries based on the socioeconomic indicators. The clustering is implemented using PCA. Along the vertical axis, we have the distance between the two clusters being merged. The clusters are merged in agglomerative manner from bottom to top. Cutting the tree at the height of 8 results in four clusters highlighted. Countries having the same color shade belong to the same cluster. Uganda is a single cluster element.}
  \label{fig: socioecon}
\end{sidewaysfigure}

Using the rule of thumb and considering the merging where the clusters have the most significant gap, we cut the hierarchy at the height of 2.5 and generated three meaningful sets as shown in Figure \ref{fig: dpca dendo conta}. The same intuition was also applied in developing the dendrogram in Figure \ref{fig: dpca cont socio} and Fig. \ref{fig: socioecon} by cutting the hierarchy at a certain level.

The first cluster consists of only one country, Uganda. Compared to other countries, it has the lowest social protection, the highest population growth, and the highest foreign direct investment, with an Internet penetration of only 20 percent. It can be seen from Fig. \ref{fig: socioecon} that it has a different economic performance. However, the country has the same contact patterns as Sierra Leone (Fig. \ref{fig: dpca dendo conta}).

The countries of the second cluster, shown in orange as shown in Figure \ref{fig: dpca cont socio}, mainly include countries from Southern Africa and a few from West and East Africa. These countries have high labor force participation (more than 0.6) and low social security coverage. The socioeconomic indicators of these countries show average performance, indicated by average life expectancy. Mozambique, Zimbabwe, and Botswana have high HIV incidence rates and slum populations. South Africa is the most developed country in this cluster, with good access to health and social security, good educational infrastructure, and low levels of food insecurity. It has a labor force participation rate of 53 percent, the lowest slum population compared to highly urbanized countries, and is the most industrialized and therefore characterized by the highest carbon emissions \cite{Z}.

The green-colored countries of the third cluster consist of Nigeria, Ethiopia, Lesotho, Liberia, and Mauritania. These countries show similar economic indicators (Fig. \ref{fig: socioecon}). For example, they have an average electricity connection between 47 and 55 percent, which has led to low internet penetration. The maternal modalities in this cluster are higher than in the rest of the groups \cite{Z}. These countries have the same social interactions (Figure \ref{fig: dpca dendo conta}).

The fourth cluster, colored in red, mainly includes countries from North Africa (Algeria, Morocco, Tunisia, and Egypt), island states (Seychelles and Mauritius), and others from the rest of Africa. Algeria, Morocco, Tunisia, and Seychelles have the highest life expectancy at 77 years, indicating a good standard of living. Most countries in this cluster show similar socioeconomic indicators (Fig. \ref{fig: socioecon}) and contact patterns (Fig. \ref{fig: dpca dendo conta}). Table \ref{tbl} provides the final list of clusters. Uganda, is omitted since it is a one-element cluster.

\begin{table}[h]
\caption{Clusters produced by the agglomerative hierarchical clustering algorithm based on social contact patterns and socioeconomic indicators for the 32 African countries applying PCA and projection $(\mathrm{2D})^{2} PCA)$. Cluster 1, 2 and 3 have 9, 5 and 17 elements, respectively. Uganda is one element cluster and thus it is omitted in this listing.}
\begin{tabular}{||p{4cm}|p{2.8cm} |p{5.5cm}||} 
\hline
{} & {} & {}\\
Cluster 1 & Cluster 2 & Cluster 3\\[1ex]
\hline
 Botswana, Mozambique, South Africa, Burkina Faso, Guinea, Benin, Rwanda, Niger, Zimbabwe. & Nigeria, Ethiopia, Lesotho, Liberia, Mauritania. & Sao Tome, Seychelles, Senegal, Tanzania, Egypt, Mauritius, Ghana, Tunisia, Sierra Leone, Kenya, Zambia, Algeria, Congo, Cape Verde, Namibia, Cameroon, Morocco.    \\[1ex]
\hline
\end{tabular}
\label{tbl}
\end{table}

\begin{figure}{}
     \centering
     \begin{subfigure}[b]{0.33\textwidth}
         \centering
         \includegraphics[width=\textwidth]{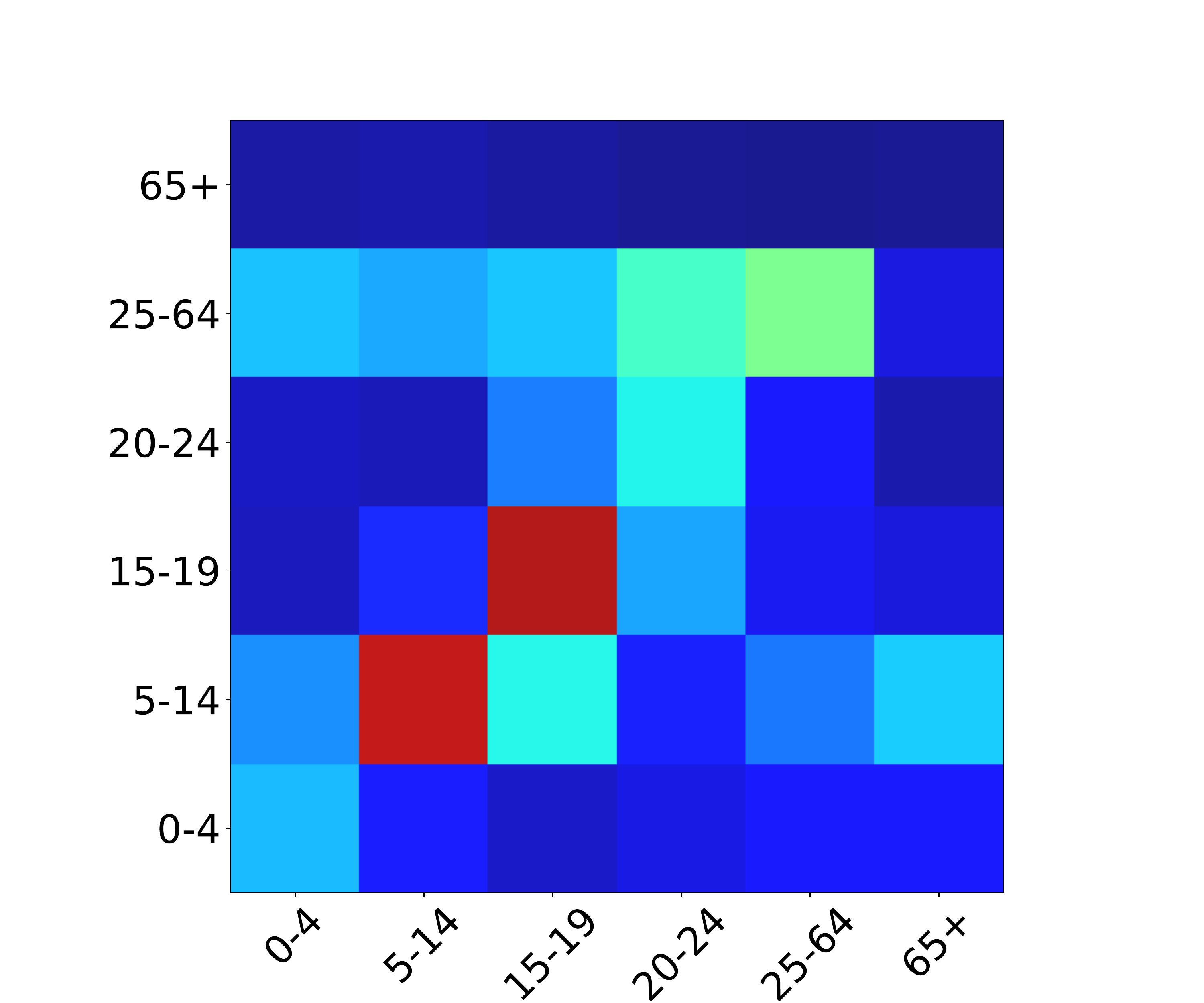}
         \caption{}
         \label{fig: cluster1}
     \end{subfigure}
     \hfill
     \begin{subfigure}[b]{0.33\textwidth}
         \centering
         \includegraphics[width=\textwidth]{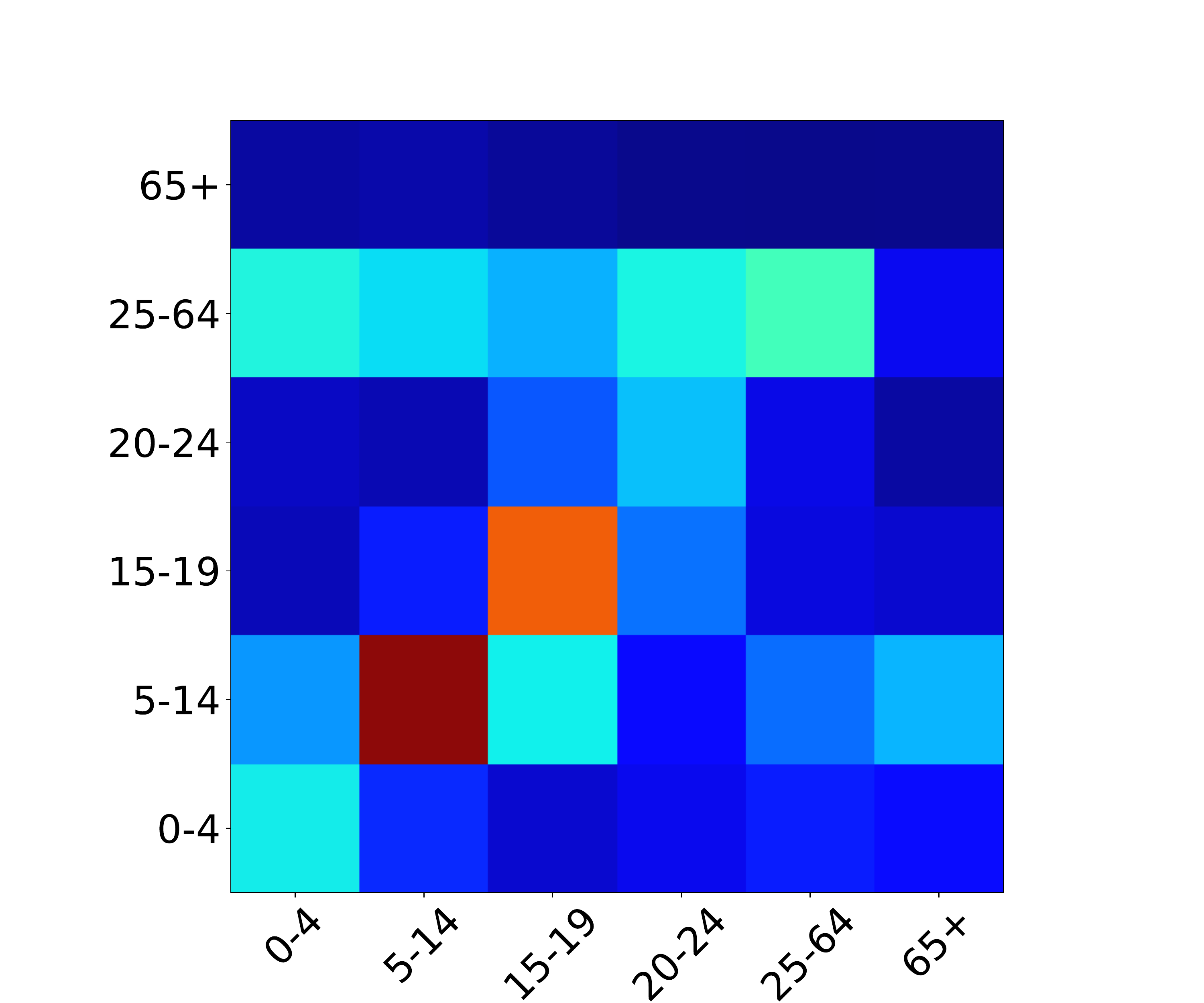}
         \caption{}
         \label{fig: cluster2}
     \end{subfigure}  
     \begin{subfigure}[b]{0.32\textwidth}
         \centering
         \includegraphics[width=\textwidth]{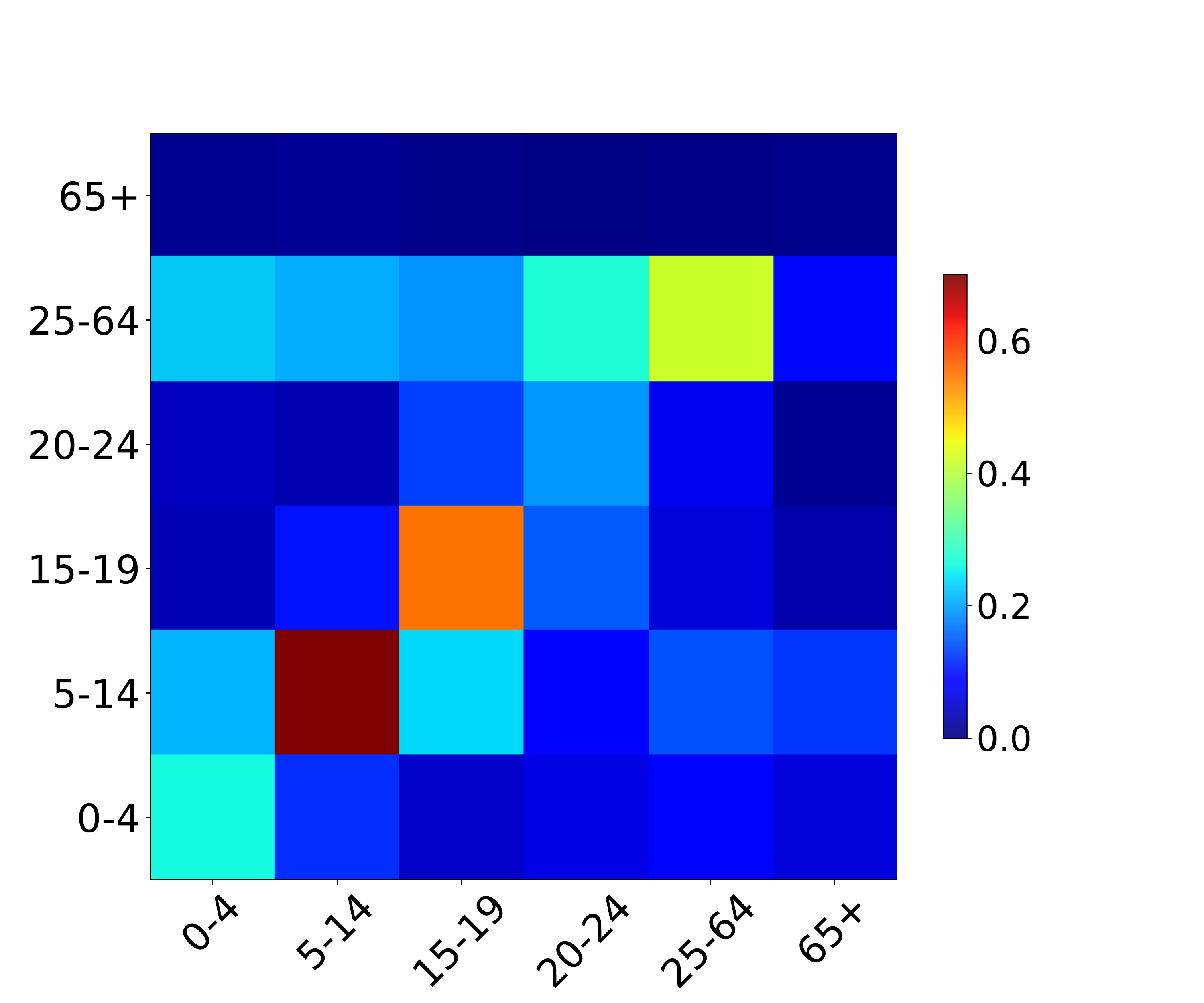}
         \caption{}
         \label{fig: cluster3}
     \end{subfigure}
         \caption{Standardized contact matrices correspond to Zimbabwe (a), Mauritania (b), and Sierra Leone (c) selected from the three clusters based on the contact matrices in Fig. \ref{fig: dpca cont socio}. These countries are located in the middle of the clusters. The pattern for the matrices are different. Cluster 1 will have the same pattern as (a), countries in cluster 2 as (b), (c) will be a representative of patterns in cluster 3.}
        \label{fig: clusters}
\end{figure}

Fig. \ref{fig: clusters} shows an element from each of the groups computed from reduced vectors by considering only the social contact patterns. These countries (Zimbabwe, Mauritania, and Sierra Leone) were chosen because they are in the middle of the blocks along the horizontal axis \ref{fig: dpca cont socio}. Looking at these matrices, we see differences between the matrices along the main diagonal and the contacts with the working-age population. In Zimbabwe, interaction is strong in primary and secondary school, with fewer intergenerational contacts. The working-age population has strong interactions with all age groups except the elderly. Mauritania and Sierra Leone show concentrated connections between the primary school but fewer contacts with the secondary school age groups. Zimbabwe reports some contacts for tertiary education, but these contacts are less for Sierra Leone. The connections with the workplace elements differ in all three samples. The lack of contact with the 65+ age group in the three samples is evidence of the low life expectancy in the region. Looking at the matrices, the countries in cluster 1 show the same patterns as Zimbabwe. The countries in cluster 2 as Mauritania and Sierra Leone, will represent the patterns of the countries in cluster 3.

\section{Discussion}

In the case of an emerging infectious disease, such as COVID-19, the transmission does not occur in all countries simultaneously. Understanding the countries to watch ahead of a pandemic is crucial. It is essential to identify those countries that show comparable patterns and consequently can develop similar NPIs. Using hierarchical clustering, we categorized the African countries into groups based on associated social contacts and socioeconomic variables. Using four principal components, we found four clusters that show high similarity in social connections and socioeconomic patterns within the cluster and lower similarity between groups. 
The country in the first cluster (Uganda) is still developing. It faces food insecurity, high unemployment rates, high population growth, weak power connections, and low Internet penetration rates. The cluster relies on manufacturing and agriculture, so most people's livelihoods depend on daily activities. Implementing an NPIs strategy such as mandatory social distancing and lockdown will disrupt supply chains and negatively impact the economy and people's livelihoods. In a pandemic situation, the government of such countries should strike a balance between the pandemic and people's living standards. In this cluster, the partial adoption of NPIs can balance the country's economy, people's livelihood, and pandemic spread. Alternatively, governments can implement NPIs but sort funds from international organizations like the World Bank to cushion vulnerable people during the pandemic.

The countries in the second cluster are mainly from southern Africa. The Human Development Index has ranked most countries in this cluster as the least developed countries. Countries have reported high rates of food insecurity, poverty, unemployment, population growth, mortality rates, low literacy rates, weak power connections, and internet penetration rates. Countries like South Africa, Mozambique, Zimbabwe, and Botswana have very high HIV incidence rates. In this cluster, there are strong contacts between children of school age at home and few contacts of working age  (Fig. \ref{fig: cluster1}).
During a pandemic, school closures and partial adoption of the NPIs strategy to balance the country's economy and livelihood can effectively control the pandemic.

In the third cluster, we have countries in Africa that are slightly developed. The countries have a moderate labor force and low unemployment rates, except for Lesotho, which reports a somewhat higher rate. Most interactions in these countries (Fig. \ref{fig: cluster2}) take place in schools (primary level). We have few contacts at work, home, and other areas like public places and communities. School closures and working from home can be efficient in these countries to contain a pandemic.

Most countries in the fourth cluster have average socioeconomic performance and high labor force participation rates. Some countries, especially the northern and island countries, have better living standards and lower poverty rates. However, northern countries have lower labor force participation rates of less than 50 percent and high unemployment rates. Most countries in this cluster have high secondary education enrollment, which can be seen from the standardized contact matrix, where we have concentrated contacts between the young children, indicating interactions in their respective schools (Fig. \ref{fig: cluster3}). In this cluster, we have strong contacts at school and a few in other settings, including the workplace. Introducing NPIs measures such as home working alongside school closures may prove effective in containing the spread of the pandemic in this cluster. 

In contrast to more classical regression techniques that try to find the relationship between features and a single target, clustering algorithms using dimensionality reduction techniques offer the possibility to visualize the country-specific social contacts and variable socioeconomic heterogeneity and to identify the clusters that show similar patterns. The results can be useful for policymakers and governments in designing country-specific NPIs to reduce the transmission of pandemics in the countries. This is crucial given the lack of public health interventions such as vaccination and prior knowledge of the countries that need monitoring and comparison.

The study has some limitations. We used estimated contact matrices of \cite{PK}, which depend on the sample size available for each country and the statistical method. More smoothing was applied in countries with small sample sizes, affecting the estimates. The quality of the estimates affects the defined distance and, thus, the clusters. Some information needs to be recovered in the dimensionality reduction techniques, but the methods preserve the patterns of the distance matrix. The 28 socioeconomic indicators on the \cite{Z} website were presented for a different year. Some indicators were reported for 2015, e.g., the Poverty rate. Literacy, research, and development expenditures were presented for 2013, along with bribery rates (2010). Some variables were unavailable in the database for some countries included in our analysis.

\section{Conclusions}

PCA enabled us to identify the critical components from the socioeconomic variables that stored 61.6 percent of information across countries. Using agglomerative hierarchical clustering, we captured complex social contact patterns associated with socioeconomic factors shared by governments to generate countries with similar practices that explain the course of an infectious disease. Our results can be helpful for policymakers to design the NPIs aimed at reducing transmission in the selected countries.

\section{Data Availability}
Code and country-level data used for this analysis are available on GitLab:
\url{https://github.com/Evanskorir/African-social-contact-patterns}

\section{Acknowledgement}
The Stipendium Hungaricum Scholarship Program with Application No. 118250 supported E.K.K.\\
 Z.V was partially supported by the National Laboratory for Health Security Program RRF-2.31-21-2022-00006 and by project TKP2021-NVA-09, implemented with the support provided by the Ministry of Innovation and Technology of Hungary from the National Research, Development and Innovation Fund, financed under the TKP2021-NVA funding scheme.

\section{Appendix}

\subsection{List of the African countries}

We considered the following African countries in the study:\\
 Algeria, South Africa, Morocco, Tunisia, Cameroon, Mauritius, Seychelles, Uganda, Zambia, Ghana, Namibia, Congo, Guinea, Botswana, Benin, Mozambique, Rwanda, Sierra Leone, Tanzania, Ethiopia, Lesotho, Sao Tome, Burkina Faso, Liberia, Nigeria, Egypt, Senegal, Kenya, Zimbabwe, Niger, Cape Verde, Mauritania.
  
\subsection{List of the socioeconomic indicators}

We considered the following indicators in the study:\\
Population Density, Urban Population, Labour Force Participation,
School enrollment, Social Protection Coverage, Poverty head count ratio, GDP Per Person Employed, Food insecurity prevalence, Life expectancy, Population Growth Rate, GDP Growth Rate, Unemployment Rate, Inflation Rate, Personal remittances, Electricity Access, Internet penetration Rate, Foreign Direct Investment, Maternal Mortality Ratio, Under-five Mortality Rate, HIV Incidence Rate, Fertility Rate, Literacy Rate, Value added Percent, Slums Proportion, Bribery incidence, Exports, Debt service, Tax revenue.

\end{document}